\documentclass{aptpub}

\usepackage[colorlinks = true,linkcolor = blue,urlcolor  = blue,citecolor = blue,anchorcolor = blue]{hyperref}
\authornames{T.~BRITTON AND D.~ZHANG} 
\shorttitle{Epidemic models with digital and manual contact tracing} 

\newcommand{\addresscurrent}[1]{\footnote{\hspace*{-14pt}$^{**}\,$Current address: #1}\par}

\usepackage{dcolumn}        
\usepackage{booktabs}       
\usepackage{paralist}       
\usepackage{graphicx} 
\usepackage{subcaption}

\begin{document}

\title{Epidemic models with digital and manual contact tracing} 

\authorone[Stockholm University]{Tom Britton} 
\authortwo[Stockholm University]{Dongni Zhang}

\addressone{Department of Mathematics, Stockholm University, 106 91 Stockholm, Sweden.} 
\addresscurrent{Department of Health, Medicine and Caring Sciences, Link{\"o}ping University, 581 83 Link{\"o}ping, Sweden.} 
\emailone{tom.britton@math.su.se} 
\emailtwo{dongni.zhang@liu.se (corresponding author)} 
\begin{abstract}
We analyze a Markovian SIR epidemic model where individuals either recover naturally or are diagnosed, leading to isolation and potential contact tracing. Our focus is on digital contact tracing via a tracing app, considering both its standalone use and combination with manual tracing. We prove that as the population size $n$ grows large, the epidemic process converges to a limiting process, which, unlike typical epidemic models, is not a branching process due to dependencies created by contact tracing. However, by grouping to-be-traced individuals into macro-individuals, we derive a multi-type branching process interpretation, allowing computation of the reproduction number $R$. This is then converted to an individual reproduction number $R^{(ind)}$, which, contrary to $R$, decays monotonically with the fraction of app-users while both share the same threshold at 1. Finally, we compare digital (only) contact tracing and manual (only) contact tracing, proving that the critical fraction app-users $\pi_c$ required for $R=1$ is higher than the critical fraction manually contact traced $p_c$ for manual tracing. 
\end{abstract}



\section{Introduction}\label{sec:intro}
One of the main reasons for modelling epidemics is to understand the effect of different preventive measures, such as vaccination and non-pharmaceutical interventions (NPI). During the recent Covid-19 pandemic, research about the effect of various NPIs, including social distancing, case isolation, contact tracing, lockdown, etc., has received much attention. Some related papers (e.g. \cite{ferguson_report_2020, flaxman_estimating_2020,longini_containing_2005}) have been highly influential on public health policies. 

This paper focuses on the preventive measure: ``\textit{contact tracing}". 
Contact tracing is traditionally performed by public health agency officers who interview diagnosed individuals and then call the named contacts to advise them to self-quarantine and to test, which is often referred to as the \textit{manual contact tracing}. This type of contact tracing has been successful in reducing transmission in many epidemics like Ebola \cite{ebola_ct}, SARS \cite{lloyd-smith_curtailing_2003}, Influenza \cite{ref1}, and Measles \cite{liu_role_2015}. 
In the present paper, we instead investigate a different type of contact tracing:  \textit{digital contact tracing}. In the case of Covid-19, \cite{ferretti_quantifying_2020} suggested that manual contact tracing is not rapid enough to control the Covid-19 epidemic, and the new idea of developing a ``\textit{contact tracing app}" was invented: once an app-user is diagnosed, a warning message will be instantaneously sent out to all the app-users who recently have been nearby for a sufficient duration with the confirmed case. Such notified contacts are then advised to test and quarantine themselves. Several contact tracing apps have been developed and released (e.g. in the UK \cite{ferretti_quantifying_2020,wymant_epidemiological_2021}, in the Netherlands \cite{app_netherlands}). 

There are several papers analyzing digital contact tracing (e.g. \cite{pollmann_impact_2020,kretzschmar_isolation_2020,kucharski_effectiveness_2020,barrat_effect_2021}) focusing on simulation-based epidemic models, and the results of \cite{rizi2022epidemic} are based on a combination of simulations and mean-field analysis of the related percolation problem. Several papers on digital contact tracing analyzing Covid-19 (e.g. \cite{cencetti_digital_2021,ferretti_quantifying_2020})  are based on the mathematical model first introduced by \cite{fraser_factors_2004} in terms of recursive equations for analyzing the timing of infectiousness and the appearance of symptoms. In particular, the results from \cite{ferretti_quantifying_2020} showed that the Covid-19 epidemic was unlikely to be contained only by manual tracing. 

Modelling contact tracing is mathematically challenging and does not easily lend itself to being analyzed using differential equations. Our focus in this paper is on the beginning of an outbreak before a massive fraction have been infected. For this reason, we apply stochastic models, which are also suitable because the outcome of contact tracing itself is random rather than deterministic. The early phase of an epidemic without preventive measures can often be approximated by a suitable branching process (e.g. \cite{ball_strong_1995}) using the underlying model assumption stipulating that once an individual has been infected, it behaves (transmits) independently of its infector. When contact tracing is introduced, this independence breaks down and other methods need to be developed. In the last decades, there has been some progress in more rigorous analyzes of epidemic models with manual contact tracing, e.g. \cite{muller_contact_2000, ball_threshold_2011,barlow_branching_2020}.

In \cite{zhang_analysing_2021}, we considered a Markovian SIR epidemic model with manual contact tracing by assuming that once individuals are diagnosed, they perform manual tracing without delay, thus reaching each contact with probability $p$. 
In the present paper, we instead analyze a model with \emph{digital} contact tracing by assuming that a fraction $\pi$ of individuals installs a contact-tracing app (and follows the advice if notified). Once an app-user is diagnosed, all its app-using contacts are notified, with the result that they isolate and test themselves without delay. For those who test positive, their app-using contacts are also notified, and so on. Furthermore, for both manual and digital contact tracing, we assume that traced individuals who have recovered naturally are also identified and are subject to contact tracing. Please note the conceptual difference between digital and manual contact tracing: manual contact tracing may happen if a contacted \emph{individual} is traced (which happens with probability $p$), whereas digital contact tracing happens if both the individuals are app-users thus creating an app-using \emph{connection} (or edge).

Using large population asymptotics, we prove that the early epidemic process with digital tracing may be approximated by a two-type branching process where type-1 ``individuals" are \textit{app-using components} (in terms of app-users) and type-2 are non-app-users. We then analyze the model having both types of contact tracing in place, i.e. reaching a fraction $p$ of all contacts using manual contact tracing and additionally all app-using contacts if the tested individual is an app-user (making up a fraction $\pi$ in the community). The limiting process in this case is shown to be a two-type branching process, where both type-1 and type-2 represent the to-be-traced components - referred to as ``macro-individuals", but type-1 starts with an app-user and type-2 with a non-app-user. 

We derive reproduction numbers for the two limiting processes. These reproduction numbers are however hard to interpret, being reproduction numbers of multi-type branching processes where some of the types are not actual individuals but rather components of individuals. The reproduction numbers hence reflect the average number of new clumps infected by such clumps during the early phase of the epidemic. These average numbers need not be the same as the average number of \emph{individuals} that one typical infected \emph{individual} infects during the beginning of the outbreak, which would define the individual basic reproduction number. The important connection between such complicated reproduction numbers and the individual reproduction numbers is that they are either both above 1 (the super-critical case), both below 1 (the sub-critical case), or both exactly 1. Because the multi type reproduction numbers are hard to interpret, we also derive corresponding individual reproduction numbers.

However, individual reproduction numbers may also be defined differently depending on who is attributed to specific infections. For example, in a household epidemic model, if all infections in a household are attributed to the person infecting the household index case, then this will give rise to one individual reproduction number. In contrast, attributing infections to the person who actually infected the individual (a more natural but also less mathematically tractable definition) gives rise to the individual basic reproduction number. As a consequence, there may also be different individual reproduction numbers for a given model (relating to the time ordering of when infections occur). Fortunately, they also share the relation of all being larger than 1, smaller than 1, or equal to 1. For a further treatment of various reproduction numbers and their interpretation, we refer to \cite{ball_trapman_2016}.
An interesting observation is that the component reproduction number for the digital-only model may in fact \textit{increase} with $\pi$, the fraction of app-users. Still, it is proven that the corresponding individual reproduction number is monotonically decreasing with $\pi$ as expected.

We also compare the digital (only) contact tracing model and compare it with the manual (only) contact tracing model of \cite{zhang_analysing_2021} to see how the critical app-using fraction $\pi_c$ relates to the critical fraction $p_c$ of contacts that are manually contact traced, having all other model parameters equal. We prove that $\pi_c>p_c$ and, as a by-product, obtain a new explicit expression for $p_c$ only given implicitly in \cite{zhang_analysing_2021}.

In Section \ref{sec:model}, we define the digital contact tracing model as well as the model having both manual and digital contact tracing. {In Section \ref{sec:mainresults}, we present the main results and provide some intuition.} Next, {we give the proofs and explain more about} approximating the initial phase of the epidemic in a large population ({Section \ref{sec:early_dct}} for digital tracing only, {Section \ref{sec:early_comb}} for both types of tracing). Then, in Section \ref{sec:numerical}, we illustrate our results numerically, compare the effect of manual and digital contact tracing and investigate their combined effect. Finally, we summarize our conclusions and discuss potential improvements of the present models in Section \ref{sec:discussion}.

\section{{Model description}}
\label{sec:model}

\subsection{The epidemic model}\label{sec:epimodel}
First, we consider a continuous-time Markovian SIR (\textbf{S}usceptible $\rightarrow$ \textbf{I}nfectious $\rightarrow$ \textbf{R}ecovered) epidemic spreading in a closed, homogeneously mixing population of fixed size $n$. Initially, there is one infectious individual, and the rest are susceptible. An infectious individual makes contact at rate $\beta$, each time with an individual chosen uniformly at random from the population. The individual-to-individual contact rate is thus ${\beta}/{n}$. If the contacted individual is susceptible, they become infected immediately; otherwise, nothing happens. Each infectious individual recovers naturally at rate $\gamma$ and plays no further part in the disease spreading (thus having an infectious period, denoted by $T_{I}$, exponentially distributed with mean ${1}/{\gamma}$).

 To incorporate contact tracing, we add a testing strategy to this SIR model. We assume that once infected, infectious individuals test positive and are immediately isolated (referred to as ``diagnosed") after random periods, which are independent and identically distributed according to a random variable $T_{D} \sim Exp(\delta)$. Therefore, infectious individuals can stop spreading the disease after a random time $min\{T_{I}, T_{D}\} \sim Exp(\gamma+\delta)$, either through natural recovery (rate $\gamma$) or diagnosis (rate $\delta$). All the random quantities described are assumed to be mutually independent. The epidemic stops when no infectious individuals remain in the population.

In this epidemic model (with testing but without contact tracing), each infectious period follows an exponential distribution with intensity $(\gamma+\delta).$ Consequently, it is straightforward to compute the mean number of secondary infections produced by a single infective before recovery or diagnosis:
\begin{equation}\label{eq:expression_R_0}
    R_0 = \frac{\beta}{\gamma + \delta}.
\end{equation}

For emphasis on the effect of contact tracing, we refer to this as the \textit{basic reproduction number} $R_0$, though more correctly, $R_0$ should represent the situation without testing (so $\delta\to 0$ and $R_0$ will be equal to ${\beta}/{\gamma}$). However, we use the expression above to study the impact of contact tracing versus no contact tracing, assuming the testing rate remains unchanged.

\subsection{The epidemic model with digital contact tracing}\label{sec:dct_model}

We add digital contact tracing to the epidemic model in Section \ref{sec:epimodel} by assuming that a fraction $\pi$ of individuals use the tracing app (and follow the recommendations) and that individuals mix uniformly irrespective of using the app or not, {i.e. each individual is an app-user with probability $\pi$ independently. In the following text, we assume that $0<\pi \le 1$ unless otherwise specified. Clearly, when $\pi=0$, there is no digital tracing, so we would have the same epidemic model as in Section \ref{sec:epimodel}.}

{For the infectious individuals (either app-users or non-app-users), they recover naturally at rate $\gamma$ and are diagnosed at rate $\delta$. The digital contact tracing procedure is defined as follows.} If and when an infectious app-user is diagnosed, all app-using contacts of the infective will be notified, with the
result that they test and isolate themselves without delay, so stop spreading if infectious. The app-using contacts who have been infected, including those who have recovered, are then also assumed to trigger digital tracing among their app-using contacts without delay according to the same procedure, and so on.

\subsection{The epidemic model with digital and manual contact tracing}\label{sec:comb_model}

Finally, we describe the epidemic model with digital contact tracing when there is also manual contact tracing in place and call this the \textit{combined model}. First, {following the assumption in Section \ref{sec:dct_model},} a community fraction $\pi$ are app-users and follow the instructions. If such an app-user is diagnosed, all of its app-using contacts (on average being a fraction $\pi$ of all contacts) will be traced.
To this, we now add manual contact tracing by assuming that once an infectious individual is diagnosed, app-user or not, each of its contacts will be reached independently by manual tracing with probability $p$. It means that for diagnosed app-users, \emph{all} app-using contacts will be traced, and each non-app-using contact will be traced {manually} with probability $p$. For diagnosed non-app-users, each contact, app-user or not, will be traced {manually} with probability $p$. Like before, all traced individuals are tested and isolated without delay if testing positive. Traced individuals who test positive, including those who have recovered, are then contact traced in the same manner without delay, and so on. The model parameters are summarized in Table \ref{tab:parameter}.

\begin{table}[tb!]
\begin{center}
\caption{Table of Model Parameters}\label{tab:parameter}  
\begin{tabular}{@{}llll@{}}
 \toprule
Parameter & Notation \\
 \midrule
size of population & $n$ \\
transmission rate & $\beta$ \\
rate of natural recovery & $\gamma$ \\
rate of diagnosis & $\delta$ \\
community fraction using the contact tracing app & $\pi$\\
(and following the advice)& \\
probability that contact is manually traced   & $p$ \\
\bottomrule
\end{tabular}
\end{center}
\end{table}

\section{{Main results}}\label{sec:mainresults}
In this paper, our main goal is to analyze the initial phase of the epidemic having digital contact tracing. We start by showing that the epidemic converges to a limiting process, asymptotically as the population size $n$ grows to infinity.

\subsection{Digital contact tracing only}
{Let $E^{(n)}_{D}(\beta,\gamma,\delta,\pi)$ denote the epidemic model with digital contact tracing only in a population of size $n$  (defined in Section \ref{sec:dct_model}). It is convenient to consider the epidemic as a two-type process where the type of infectives corresponds to whether they are app-users or not. The limiting process, denoted $E_{D}(\beta,\gamma,\delta,\pi)$, is described as follows. Each alive (=infectious) individual gives birth to (infects) app-users at rate $\beta\pi$, to non-app-users at rate $\beta(1-\pi)$, recovers naturally at rate $\gamma$ and is diagnosed at rate $\delta$ (the individual is considered dead, not giving any more births, in both cases). Moreover, if a diagnosed individual is an app-user, all of its app-using offspring as well as its parent (in the case of an app-user) are diagnosed at the same time. Such diagnosed individuals will, in turn, lead to all app-users among its offspring and parent being diagnosed without delay, and so on.} 

The limit process in terms of app-users and non-app-users will \textit{not} be a branching process when contact tracing is in place, since an infectious app-user maybe contact traced by its app-using infector or one of its app-using infectees, thus giving rise to dependence between individuals which violates the basic assumption of independence in branching processes. But if we instead consider non-app-users as one type of individual and  ``app-using components'' (consisting of the whole transmission chain of app-users) as a second type of ``macro-individual" this process will evolve according to a two-type branching process (see Figure \ref{fig:two_type_br_digital} of the next section for an illustration). This holds because there will never be contact tracing between any combination of such individuals, thus making individuals behave independently.

Next, we study the birth and death of app-using components. Such a component starts with one infectious app-user infected by a non-app user (we call this the \textit{root} of the component). This component can then grow as the root infects other app-users, who in turn infect more app-users, and so on. Eventually, the whole app-using component will die out (i.e. stop spreading the infection) either when someone within the component is diagnosed (leading to immediate tracing and diagnosis of the whole component) or when all individuals in the component recover naturally. While the component has infectious individuals, it can also infect non-app-users. These infected non-app-users can subsequently infect other non-app-users and app-users; in the latter case, new app-using components are created.

Following  multi-type branching process theory \cite{ref2,becker_multi_type_1990} we define
\begin{equation}\label{eq:mean_matrix}
    M=\begin{pmatrix}
m_{11} & m_{12}\\
m_{21}  & m_{22} 
\end{pmatrix}
\end{equation}
as the next generation matrix of the limiting two-type branching process, where $m_{ij}$ is the expected number of secondary infections of type $j$ produced by a single infected individual of type $i$, for $i,j=1,2$ (type 1 is app-using components and type 2 non-app-users). Then, a major outbreak occurs with positive probability if and only if the largest eigenvalue of $M$ denoted $R_{D}$, is above 1 (the elements of $M$ are derived in Section \ref{sec:R_D}). We call $R_{D}$ the \textit{effective reproduction number for the epidemic with digital tracing}. 

We summarize our findings in the following proposition and corollary. 
\begin{prop}\label{prop:earlyapprox_dct}
Consider the sequence of epidemic processes with digital tracing $E^{(n)}_{D}(\beta,\gamma,\delta,\pi)$, for each population size $n $, starting with one random initial infective and the rest susceptible. Then on any finite time interval, $E^{(n)}_{D}(\beta,\gamma,\delta,\pi)$ converges in distribution to the limit process $E^{}_{D}(\beta,\gamma,\delta,\pi)$ as $n\to\infty.$ 
\end{prop}

\begin{cor}\label{col:digital}   
The limit process $E^{}_{D}(\beta,\gamma,\delta,\pi)$ is super-critical, meaning that it will grow beyond all limits with positive probability if and only if $R_D>1$, where 
\begin{equation}\label{eq:R_D}
    R_{D} = \frac{R_{0}(1-\pi)}{2}+\sqrt{\frac{( R_{0}(1-\pi))^2}{4}+f(\pi)R_{0}\pi}
\end{equation}
with 
\begin{equation}\label{eq:expression_f}
   f(\pi)= \frac{(1-\pi)}{2 \pi \delta} (\beta \pi-\gamma-\delta+\sqrt{(\beta \pi+\gamma+\delta)^2-4\beta \pi \gamma}).
\end{equation}
Let $Z^{(n)}_{D}$ be the final number of infected in the whole epidemic $E^{(n)}_{D}(\beta,\gamma,\delta,\pi)$. Then if $R_{D} > 1$, $Z^{(n)}_{D} {\to} \infty$ with positive probability, whereas if $R_{D} \leq 1$ then the final fraction getting infected $\bar Z^{(n)}_{D} = Z^{(n)}_{D}/n \overset{p}{\to} 0$ (a minor outbreak for sure).

\end{cor}

Our next result compares the effect of digital contact tracing with the effect of manual contact tracing. More precisely, we compare the current model with the manual contact model of \cite{zhang_analysing_2021}. That model has the same transmission model, no digital contact tracing, but instead a manual contact tracing in that each contact is traced upon diagnosis, independently, with probability $p$. Keeping all epidemic parameters being fixed and equal in the two models, we compare $\pi_c$ being the app-using fraction for which $R_D=1$ with $p_c$ being the critical fraction manually contact traced such that the corresponding manual reproduction number satisfies $R_M=1$. The following theorem shows that a larger app-using fraction is needed, compared to the manual tracing probability to reduce the reproduction number to 1.

\begin{thm}\label{thm:p_c_pi_c}
Let $\beta,\delta,\gamma >0$ be fixed parameters and assume $R_{0}=\beta/(\delta+\gamma)>1$ (otherwise contact tracing is unnecessary). 
  Then we have 
     $$\pi_c > p_c = \frac{\beta-(\gamma+\delta)}{\beta-\gamma}.$$
\end{thm}

The reproduction number $R_D$ cannot be interpreted as some mean number of infections caused by one individual since it is based on a two-type branching process in which one of the types is a to-be-traced component. A notable consequence of this is seen in Figure \ref{fig:RD} where we observe that $R_{D}$ is not monotonically decreasing in $\pi$: in some situations, increasing the fraction of app-users can lead to a \emph{larger} reproduction number!

In Section \ref{sec:Rind_D} we therefore 
derive an individual reproduction number $R^{(ind)}_{D}$ for the epidemic with digital tracing:
\begin{equation*}
     R^{(ind)}_{D}= \frac{\pi f(\pi)}{1-\pi +\pi f(\pi)} + (1-\pi) \frac{\beta}{\delta+\gamma}
\end{equation*}
with $f(\pi)$ given by \eqref{eq:expression_f}.
For this individual reproduction number $R^{(ind)}_{D}$, we have the following theorem showing that it is monotonically decreasing in $\pi$. The subsequent theorem shows that $R^{(ind)}_{D}$ is indeed an epidemic threshold. 

\begin{thm}\label{thm:R_D_ind_monotone}
Let $\beta,\gamma$ and $\delta$ be fixed.
Then $R^{(ind)}_{D}$ is monotonically decreasing in $\pi.$ That is, for any $0\leq \pi_{1}<\pi_{2}\leq 1$, we have $R^{(ind)}_{D}(\pi_{1}) > R^{(ind)}_{D}(\pi_{2})$.
\end{thm}

\begin{thm}\label{thm:R_D_ind}
    {Let $\beta,\gamma$ and $\delta$ be fixed. 
Then $R_D=R^{(ind)}_{D}=1$ at $\pi = \pi_c$, and $R_D>1$ as well as $R^{(ind)}_{D}>1$ if $\pi<\pi_c$ and the opposite if $\pi>\pi_c$.
}
\end{thm}
\subsection{Combined model with manual and digital contact tracing}
Finally, we state related results as above, but now for the beginning of the epidemic in the combined model of Section \ref{sec:comb_model} having digital as well as manual contact tracing, which we denote by $E^{(n)}_{DM}(\beta,\gamma,\delta,p,\pi)$. More specifically, the result proves that the combined tracing model is well approximated by the multi-type process defined in the text, and that this process when characterized by its to-be-traced components allows for a branching process description.

\begin{prop}\label{prop:earlyapprox_comb}
Consider a sequence of epidemic processes with both manual and digital tracing $E^{(n)}_{DM}(\beta,\gamma,\delta,p,\pi)$, for population size $n$, and starting with one initial infective. Then on any finite time interval, $E^{(n)}_{DM}(\beta,\gamma,\delta,p,\pi)$ converges in distribution to the limit process $E^{}_{DM}(\beta,\gamma,\delta,p,\pi)$ as $n\to\infty.$
\end{prop}
\begin{cor}\label{col:combined}
Let $Z^{(n)}_{DM}$ be the final number infected in the whole epidemic $E^{(n)}_{DM}(\beta,\gamma,\delta,p,\pi)$. Then if the reproduction number $R_{DM} > 1,$ $Z^{(n)}_{DM} {\to} \infty$ with non-zero probability; and if $R_{DM} \leq 1,$ the final fraction $\bar Z^{(n)}_{DM} = Z^{(n)}_{DM}/n \overset{p}{\to} 0$, namely there will be a minor outbreak for sure.
\end{cor}

Table \ref{tab:key_quantities} lists the above important quantities.
\begin{table}[t]
\begin{center}
\caption{Key quantities related to the epidemic processes and their corresponding limit branching processes}
\label{tab:key_quantities} 
\renewcommand\arraystretch{1.2}
\begin{tabular}{@{}ll@{}}
\toprule
Notation  &  Description\\
\midrule
$E^{(n)}_{D}(\beta,\gamma,\delta,\pi)$ & epidemic process with digital tracing\\
&(population size $n$, starting with one initial infective)\\
$E^{}_{D}(\beta,\gamma,\delta,\pi)$ & limit branching process for epidemic with digital tracing\\ 
$R_{D}$ & reproduction number for the limit process $E^{}_{D}(\beta,\gamma,\delta,\pi)$\\
$R^{(ind)}_{D}$ & individual reproduction number for the epidemic \\
&with digital tracing\\
$E^{(n)}_{DM}(\beta,\gamma,\delta,p,\pi)$ & epidemic process with digital and manual tracing\\
&(population size $n$, starting with one initial infective)\\
$E^{}_{DM}(\beta,\gamma,\delta,p,\pi)$ & limit branching process for epidemic with digital \\
&and manual tracing (combined model)\\
$R_{DM}$ & reproduction number for the limit process $E^{}_{DM}(\beta,\gamma,\delta,p,\pi)$\\
\bottomrule
\end{tabular}
\end{center}
\end{table}
\section{Approximation of the initial phase in a large community with digital tracing only}\label{sec:early_dct}

In the following subsections, we first prove Proposition \ref{prop:earlyapprox_dct} in Section \ref{sec:proof_thm_dct}, then derive the effective reproduction number for digital tracing and prove Corollary \ref{col:digital} in Section \ref{sec:R_D}. In Section \ref{sec:proof_thm_compare}, the proof of Theorem \ref{thm:p_c_pi_c} is presented, where we compare the effectiveness of manual and digital tracing analytically. Finally in Section \ref{sec:Rind_D}, we derive an individual reproduction number for digital tracing and prove Theorem \ref{thm:R_D_ind_monotone}, \ref{thm:R_D_ind}.

\subsection{{Proof of Proposition \ref{prop:earlyapprox_dct}}}\label{sec:proof_thm_dct}
We recall that for the size of population $n$, the epidemic process with digital contact tracing is denoted by $E^{(n)}_{D}(\beta,\gamma,\delta,\pi)$, where there is one initial infective. The limit process $E_{D}(\beta,\gamma,\delta,\pi)$ is defined as follows. There are two types of individuals: non- and app-users. At time $t=0$, there is one initial ancestor (who has the same type as the initial case in $E^{(n)}_{D}(\beta,\gamma,\delta,\pi)$). During the lifetime, each individual gives birth to app-users at rate $\beta\pi$ and non-app-users at rate $\beta(1-\pi)$. Individuals recover naturally at rate $\gamma$ and are diagnosed at rate $\delta.$ Once an app-user is diagnosed, all of its app-using descendants and its parent (in the case of an app-user) are diagnosed at the same time, and so on. In particular, if traced individuals have recovered, they will still be diagnosed, implying that their app-using contacts will also be traced, and so on.

{First, we note that there is a new infection, recovery, or diagnosis in the epidemic whenever a new birth, death, or removal occurs in $E_{D}(\beta,\gamma,\delta,\pi)$, respectively. The rates of recovery ($\gamma$) and diagnosis ($\delta$) are the same as those of death and removal. An infection is ``successful" only when the contacted person is susceptible. Let $S_{n}(t)$ be the number of susceptibles in the $n^{th}$ epidemic $E^{(n)}_{D}(\beta,\gamma,\delta,\pi)$ at time $t$, the probability that a given infective infects a new app-user is $\pi S_{n}(t)/n$. This implies that the rate of new infection to app-users is $\beta\pi S_{n}(t)/n$. On the other hand, the rate of birth to app-users is $\beta\pi.$ When $n$ is large and at the beginning of the epidemic, we have $S_{n}(t)\approx n,$ so the rate of infection to app-users is close to the app-using birth rate in the limit process $E_{D}(\beta,\gamma,\delta,\pi).$ Similarly, the rate of infection to non-app-users $\beta(1-\pi) S_{n}(t)/n$ is close to the rate of birth to non-app-users $\beta(1-\pi).$ Using the standard coupling arguments for epidemic models (see \cite{ball_strong_1995}), it can be shown that on any fixed time horizon the epidemic $E^{(n)}_{D}(\beta,\gamma,\delta,\pi)$ converges in distribution to the limit process $E_{D}(\beta,\gamma,\delta,\pi)$ as $n \to \infty.$}

As discussed in Section \ref{sec:mainresults}, instead of describing $E_{D}(\beta,\gamma,\delta,\pi)$ in terms of the actual individuals, we consider it as a process of {app-using components} and non-app-users, which is a two-type branching process. If there are currently $k$ alive (=infectious) app-users in an app-using component, each such individual gives birth to (infects) new app-users at rate $\beta \pi$, so the total birth rate is hence  $k\beta \pi$. Each individual recovers naturally at rate $\gamma$, so the overall recovery rate equals $k\gamma$. Finally, the whole component is diagnosed once one of the $k$ infectious individuals is diagnosed, so the total rate of diagnosis is $k \delta$. While having $k$ infectious app-users, the component generates new (infects) non-app-users at rate $k\beta(1- \pi)$, whereas the component cannot give birth to new app-using components because any infected app-user will belong to the same component. On the other hand, the non-app-users give birth to new app-using components at rate $\beta \pi$ and non-app-users at rate $\beta(1- \pi)$.

\begin{figure}[tb!]
    \begin{center}
        \includegraphics[width=80mm]{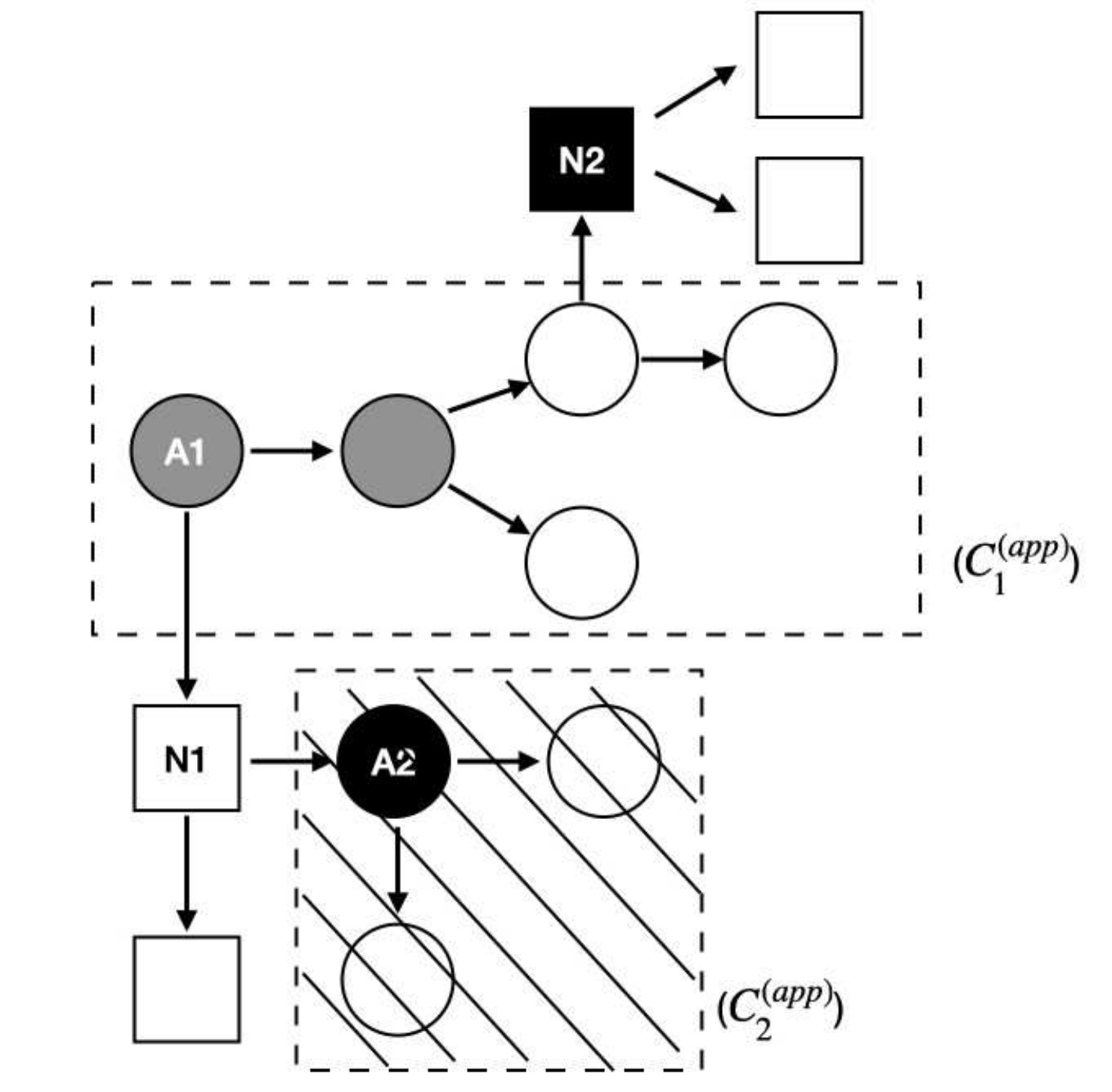}
    \caption{Example of the dynamics of app-using components and non-app-users for the model with digital contact tracing: the circle nodes are app-users and the squares non-app-users. The white, black, and gray colors correspond to ``infectious", ``diagnosed," and ``naturally recovered", respectively. The rectangular region surrounded by the dashed line symbolizes the app-using component, and the area is filled with diagonal lines when the whole component is diagnosed and removed.}
    \label{fig:two_type_br_digital}
    \end{center}
\end{figure}

In Figure \ref{fig:two_type_br_digital}, we illustrate how the two types of ``individuals": app-using components (surrounded by dashed lines) and non-app-users (square-shaped nodes) grow and eventually die out. We set the app-user $A1$ as the initial case. While the app-using component $C^{(app)}_{1}$ evolves (new app-users getting infected or existing ones recovering), the component  $C^{(app)}_{1}$ infects two non-app-users $N1$ and $N2.$ The non-app-user $N1$ infects one non-app-user and one new app-using component $C^{(app)}_{2}$ with root $A2$, whereas $N2$ infects two other non-app-users. Once the app-user $A2$ is diagnosed, the whole app-using component $C^{(app)}_{2}$ will die out.

The advantage of the approximation outlined above is that multi-type branching processes are well-studied (e.g. \cite{ref2,becker_multi_type_1990}). As a consequence, in the following section, we derive an effective reproduction number which has a threshold value of 1, determining whether a major outbreak is possible or not. 

\subsection{Proof of Corollary \ref{col:digital}}\label{sec:R_D}

From standard multi-type branching process theory \cite{ref2,becker_multi_type_1990}, the limit process $E^{}_{D}(\beta,\gamma,\delta,\pi)$ is known to be sub-critical if $R_D<1$, critical if $R_D=1$ and super-critical if $R_D > 1$, where $R_{D}$ is the largest eigenvalue of the mean-offspring matrix $M$ defined by \eqref{eq:mean_matrix}, which is given by
\begin{equation}\label{eq:eigenvalue_M}
        R_{D} = \frac{m_{11}+m_{22}}{2}+\sqrt{\frac{(m_{11}+m_{22})^2}{4}-m_{11}m_{22}+m_{12}m_{21}}.
\end{equation}
Let $Z_D$ be the total progeny of $E^{}_{D}(\beta,\gamma,\delta,\pi)$, then $Z_D < \infty$ (the process goes extinct) with probability 1 if $R_D<1$. And if $R_D>1$, $Z_D = \infty$ (the process grows beyond all limits) with a strictly positive probability. Then it follows from Proposition \ref{prop:earlyapprox_dct} that if $R_D > 1$, $Z^{(n)}_{D} \to \infty$ with non-zero probability (there may be a large outbreak of the epidemic); if $R_D \leq 1,$ $Z^{(n)}_{D} \to Z_D$, implying that the final fraction getting infected $Z^{(n)}_{D}/n \to 0$. Hence, we have proved that $R_{D}$ is an epidemic threshold. We denote the effective reproduction number for the epidemic model with digital tracing by $R_{D}$.

To get an explicit expression of $R_{D}$, we derive the elements $m_{ij}$ as follows. First, we note that the app-using components (type 1) cannot give birth to new app-using components since app-users infected by an app-using component will belong to the same component. Consequently, we have that
\begin{equation}
\label{eq:expression_m11}
    m_{11}= 0.
\end{equation}

As for non-app users, they will never be reached by digital tracing and thus spread the infection like in the epidemic without any contact tracing. That is, each infectious non-app-user infects non-app-users at rate $\beta(1-\pi)$ and app-users at rate $\beta\pi$, where the infectious periods of non-app-users are independent and identically exponentially distributed with intensity ($\delta+\gamma$). Hence, the number of app-using components (type 1) and non-app-users (type 2) infected by one typical non-app-user is geometrically distributed with support $\{0,1,2,...\}$ and mean 
\begin{equation}\label{eq:expression_m21}
    m_{21}=\frac{\beta\pi}{\delta+\gamma}=R_{0}\pi,
\end{equation}
and  
\begin{equation}\label{eq:expression_m22}
    m_{22}=\frac{\beta(1-\pi)}{\delta+\gamma}=R_{0}(1-\pi)
\end{equation}
respectively, where $R_{0}$ is given by \eqref{eq:expression_R_0}.

Finally, it remains to compute $m_{12}.$ 
Let $X(t)$ be the number of infectious app-users in an app-using component at time $t$, with $X(0)=1.$ We know that $\{X(t), t\geq 0\}$ behaves as a linear birth-death process with birth rate $\beta \pi$ (new app-users infected) and death rate  $\gamma$ (natural recovery), and this process stops when a diagnosis event happens at rate $\delta$ (the whole component is then diagnosed). Further, at time $s$, the component gives birth to new components at rate $X(s)\beta(1-\pi).$ 
Let $T$ be the random variable representing the time of first diagnosis. Then 
\begin{equation*}
    \begin{split}
        m_{12} & =  \mathbb{E}\left[\int_{0}^{T} \beta(1-\pi) X(t) dt\right]\\
        & = \beta(1-\pi) \mathbb{E}[A],
    \end{split}
\end{equation*}
where we set 
$$A := \int_{0}^{T} X(t) dt.$$

The process $\{X(t)\}$ is piecewise constant, and the jump rate of $X(t)$ is proportional to $X(t)$. More precisely, the contribution to the integral between two jumps $T_i$ and $T_{i+1}$ (prior to $T$) is $X(T_i )Y$, where $Y\sim Exp((\beta \pi + \gamma + \delta)X(T_i))$. But this contribution is simply $1\cdot \tilde Y$ where  $\tilde Y\sim Exp(\beta \pi + \gamma + \delta)$.
For this reason, $A$ equals the integral of 1 integrated up until the first time a new process $\{\Tilde{X}(t), t\geq 0\}$, which is a random \textit{time-change} version of the process $\{X(t),t\geq 0\}$, reaches the value 0. More precisely, the time span of $\{\Tilde{X}(t)\}$ is slowed down by $\{X(t)^{-1}\}$ at time $t$. This new process $\{\Tilde{X}(t), t\geq 0\}$ is hence a continuous-time random walk going up at $\beta \pi$, going down at $\gamma$ and absorbing at zero. Let $\Tilde{T}_{rw}=\inf\{t>0: \Tilde{X}(t)=0\}$, and 
$\Tilde{T}=\Tilde{T}_{rw} \land \Tilde{T}_{D},$ where $\Tilde{T}_{D} \sim Exp(\delta)$ which is independent of $\{\Tilde{X}(t), t\geq 0\}$ and hence independent of $\Tilde{T}_{rw}$. Then it follows that $A=\int_0^T X(t)dt$ and $\tilde T=\int_0^{\tilde T} 1dt$ have the same distribution. As a consequence, 
\begin{equation*}
\begin{split}
     \mathbb{E}[A] &=  \mathbb{E}[\Tilde{T}] \\
     &= \int_{0}^{\infty} \mathbb{P}(\Tilde{T}_{rw} >t) \mathbb{P}(\Tilde{T}_{D} >t)dt \\
     &= \int_{0}^{\infty} \mathbb{P}(\Tilde{T}_{rw} >t) e^{-\delta t}dt \\
     &= \frac{1-\phi(\delta)}{\delta},
\end{split}
\end{equation*}
where $\phi(\delta) := \mathbb{E}[e^{-\delta \Tilde{T}_{rw}}]$.

{Further, if we condition on the first jump of $\{\Tilde{X}(t), t\geq 0\}$, we get the following equation of $\phi(\delta)$:
\begin{equation*}
    \phi(\delta) = \frac{\beta \pi+\gamma}{\beta \pi+\gamma+\delta} (\frac{\beta \pi}{\beta \pi+\gamma}\phi(\delta)^2+\frac{\gamma}{\beta \pi+\gamma})
\end{equation*}
which has two roots $\phi(\delta)=(\beta \pi+\gamma+\delta\pm\sqrt{(\beta \pi+\gamma+\delta)^2-4\beta \pi \gamma})/{2\beta \pi}$. Since  $\phi(0)=1$, we get
\begin{equation*}
    \phi(\delta) = \frac{\beta \pi+\gamma+\delta-\sqrt{(\beta \pi+\gamma+\delta)^2-4\beta \pi \gamma}}{2\beta \pi},
\end{equation*}
implying that
\begin{equation*}
     \mathbb{E}[A] = \frac{\beta \pi-\gamma-\delta+\sqrt{(\beta \pi+\gamma+\delta)^2-4\beta \pi \gamma}}{2\beta \pi \delta}.
\end{equation*}
Hence we obtain
\begin{equation*}
    m_{12} = f(\pi):= \beta(1-\pi) \mathbb{E}[A]
\end{equation*}
for $0 < \pi \leq 1$.}

In conclusion, the effective component reproduction number $R_{D}$ for the model with digital tracing only is the largest eigenvalue to the matrix $M=(m_{ij})$, which is given by \eqref{eq:R_D}.
\begin{rem}\label{rem:R_D_continuous}
{It can be easily checked that 
\begin{align*}
    \lim_{\pi \to 0} f(\pi) = \frac{\beta(\gamma+\delta)+\beta\delta-\beta\gamma}{2\delta (\gamma+\delta)}=\frac{\beta}{\gamma+\delta} = R_{0}. 
\end{align*}
This implies that the reproduction number $R_{D}$ converges to the basic reproduction number $R_{0}$ when $\pi\to 0$ as expected: $$\lim_{\pi \to 0} R_{D}(\pi)=R_{0}=\frac{\beta}{\gamma+\delta}.$$ 
}
\end{rem}
\begin{rem}\label{rem:R_M}
The analytic expression for $f(\pi)$ can be used to obtain an analytic expression also for the manual tracing model considered in \cite{zhang_analysing_2021}. A related birth and death process was studied there. Considering such a to-be-traced component with size $k$, the component would increase by 1 at rate $k\beta p$, decrease by 1 at rate $k\gamma$, and jump to 0 at rate $k\delta.$ And before jumping, the component gave birth to new components with rate $k\beta (1-p).$ Using similar arguments as above, we hence get the following explicit expression for the reproduction number (which in \cite{zhang_analysing_2021} was expressed using an infinite sum):
\begin{equation}\label{eq:expression_R_M}
R^{(c)}_{M}(p) = f(\pi=p)= \frac{1-p}{2 p \delta} (\beta p-\gamma-\delta+\sqrt{(\beta p+\gamma+\delta)^2-4\beta p \gamma}).
\end{equation}
\end{rem}

\subsection{Proof of Theorem \ref{thm:p_c_pi_c}}\label{sec:proof_thm_compare}
In this section, we compare the separate effects of manual and digital contact
tracing by proving Theorem \ref{thm:p_c_pi_c}. We fix the epidemic parameters $\beta,\delta, \gamma$ such that $\beta/(\delta+\gamma)>1$, and consider the reproduction number for manual tracing $R^{(c)}_{M}(p)$ and the one for digital tracing $R^{}_{D}(\pi)$ as function of $p$ and $\pi$ respectively. We first state and prove two lemmas.

\begin{lem}\label{lem:1}
Suppose $R_0=\beta/(\delta+\gamma)>1$. Then the critical manual reporting probability $p_c$ (satisfying $R^{(c)}_{M}(p_c)=1$) is given by 
\begin{equation}
    p_{c} = \frac{\beta-(\gamma+\delta)}{\beta-\gamma}. 
\end{equation}
\end{lem}
\begin{proof}[Proof]

First, we show the existence of the critical probability $p_c$. Based on Remark \ref{rem:R_D_continuous} and \ref{rem:R_M}, $R^{(c)}_{M}(p)$ is continuous in $p$ with  $R^{(c)}_{M}(0)=R_{0}>1$ and $R^{(c)}_{M}(1)=0<1$. By the intermediate value theorem, there exists at least one $p_c$ in $(0,1)$ such that $R^{(c)}_{M}(p_c)=1.$

Next, we actually compute $p_c$ and show there is only one $p_c$. Without loss of generality, we can assume that $\gamma+\delta=1$, then the formula of $R^{(c)}_{M} $ in (\ref{eq:expression_R_M}) becomes
\begin{align*}
    R^{(c)}_{M}(p) &=  \frac{(1-p)}{2 p (1-\gamma)} (\beta p-1+\sqrt{(\beta p+1)^2-4\beta p \gamma})\\
&=  \frac{1}{2(1-\gamma)} \frac{(1-p)}{p } (\beta p-1+\sqrt{(\beta p-1)^2+4\beta p (1-\gamma)})\\
&= \frac{(1-p)}{2(1-\gamma)}(\beta-\frac{1}{p}+\sqrt{(\beta-\frac{1}{p})^2+4\beta (1-\gamma)\frac{1}{p}}).
\end{align*}
It follows by standard calculations that $R^{(c)}_{M}(p_{c}) = 1$ is equivalent to 
\begin{equation*}
       p_{c} = \frac{\beta -1}{\beta-\gamma}. 
\end{equation*}
We recall that the parameters have been scaled so far. Back to the original scale, we have
\begin{equation*}
    p_{c} = \frac{\beta -(\gamma+\delta)}{\beta-\gamma}.
\end{equation*}

\end{proof}

\begin{lem}\label{lem:3}
{Let $p_c$ be the critical manual reporting probability such that $R^{(c)}_{M}(p_c)=1.$ Then there exists at least one $\pi_c$ lying in $(p_c, 1)$ that satisfies $R^{}_{D}(\pi_c)=1.$}
\end{lem}
\begin{proof}[Proof]
First, we know that $R^{}_{D}(\pi)$ is continuous in $\pi$ and when $\pi=1$, $R^{}_{D}=0.$ Letting $\pi = p_c$ in (\ref{eq:R_D}) gives 
\begin{equation}\label{eq1}
    R_{D}(p_c) = \frac{R_{0}(1-p_c)}{2}+\sqrt{\frac{( R_{0}(1-p_c))^2}{4}+R_{0}p_c} > 1
\end{equation}
(the first equality follows from that $f(p_c)= R^{(c)}_{M}(p_c)=1$).
Then it follows by intermediate value theorem that there exists at least one $\pi_c$ lying in $(p_c, 1)$ that satisfies $R^{}_{D}(\pi_c)=1.$
\end{proof}

\begin{lem}\label{lem:4}
{For any critical app-using fraction $\pi_c$ such that $R^{}_{D}(\pi_c)=1,$ the reproduction number with manual tracing at $p=\pi_c$ is below 1:
\begin{equation*}
    R^{(c)}_{M}(\pi_c) < 1. 
\end{equation*}}
\end{lem}
\begin{proof}[Proof]
{$R_{D}(\pi_c)=1$ is equivalent to
    \begin{equation*}
        R^{(c)}_{M}(\pi_c)= \frac{1-R_{0}(1-\pi_{c})}{R_{0}\pi_{c}}=1-\frac{R_{0}-1}{R_{0}\pi_{c}}.
    \end{equation*}
    On the other hand, we get $({R_{0}-1})/{R_{0}\pi_{c}} >0$, since $R_{0}>1$ by assumption. This implies that $R^{(c)}_{M}(\pi_c)$ is smaller than 1.}
\end{proof}

{We are now ready to prove Theorem \ref{thm:p_c_pi_c}.}

\begin{proof}[Proof of Theorem \ref{thm:p_c_pi_c}]
{By Lemma \ref{lem:3}, we know that there exists at least one $\pi_{c}$ lying between $p_{c}$ and 1. Suppose there exists a $\pi'$ smaller than $p_c$ ($\pi' <p_{c}$) and $R_{D}(\pi =\pi') =1. $ Then by Lemma \ref{lem:4}, the reproduction number with manual tracing at point $p=\pi'$ 
$$R^{(c)}_{M}(\pi') < 1.$$}

{Again using the fact that $R^{(c)}_{M}(p=0)=R_{0}>1$ and $R^{(c)}_{M}(p)$ is continuous in $p,$ by intermediate value theorem, there shall exist one $p'$ lying in $(0,\pi')$ such that $$R^{(c)}_{M}(p') = 1.$$ On the other hand, by Lemma \ref{lem:1}, we have shown that $p'=p_{c},$ implying that $\pi' > p_{c}$ which leads to a contradiction. Hence, there does not exist such $\pi'$ which is smaller than $p_{c}$ and $R_{D}(\pi') = 1.$ }
    
\end{proof}
\subsection{The effective individual reproduction number \texorpdfstring{$R^{(ind)}_{D}$}{}}\label{sec:Rind_D}
The interpretation of $R_{D}$ is complicated in that the two types in the limiting branching process are different types of objects: type-2 individuals are actual individuals (non-app-users), but type-1 are ``macro individuals": the to-be-traced app-using components. In this section, our aim is therefore to convert $R_{D}$ to an individual reproduction number {$R^{(ind)}_{D}$}. 

 We first derive $R^{(ind)}_{D-A}$, the average number of individuals (of any type)  infected by a newly infected app-user, and $R^{(ind)}_{D-N}$ the average number of individuals (of any type)  infected by a newly infected non-app-user. We note that each infection, irrespective of the type of infector, is with an app-user with probability $\pi$ and a non-app-user with probability ($1-\pi$). This implies that the $R^{(ind)}_{D}$ can be expressed as the sum:
\begin{equation}\label{eq:expression_R_ind_pi}
    R^{(ind)}_{D} = \pi R^{(ind)}_{D-A} + (1-\pi) R^{(ind)}_{D-N} .
\end{equation}
It remains to compute $R^{(ind)}_{D-A}$ and $R^{(ind)}_{D-N}$. From the previous results in Section \ref{sec:R_D}, we see that a given app-using component infects, on average, $m_{12}$ number of non-app-users. Let $\mu^{(A)}_{c}$ be the mean number of app-users that ever get infected in an app-using component. All but one of these $\mu^{(A)}_{c}$ are infected within the component, and there are, on average, $m_{12}$ number of ``external" infections. So in total, there are on average $(\mu^{(A)}_{c} -1)+m_{12}$ number of infections by an app-using component. As a consequence, the average number of infections \emph{per app-user}, $R^{(ind)}_{D-A}$, equals 
\begin{equation}\label{eq:expression_R_App_1}
    R^{(ind)}_{D-A} = \frac{(\mu^{(A)}_{c} -1)+m_{12}}{\mu^{(A)}_{c}},
\end{equation}
where we recall from (\ref{eq:expression_f}) that 
$m_{12} = f(\pi).$

\begin{figure}[tb!]
    \begin{center}
        \includegraphics[width=.9\textwidth]{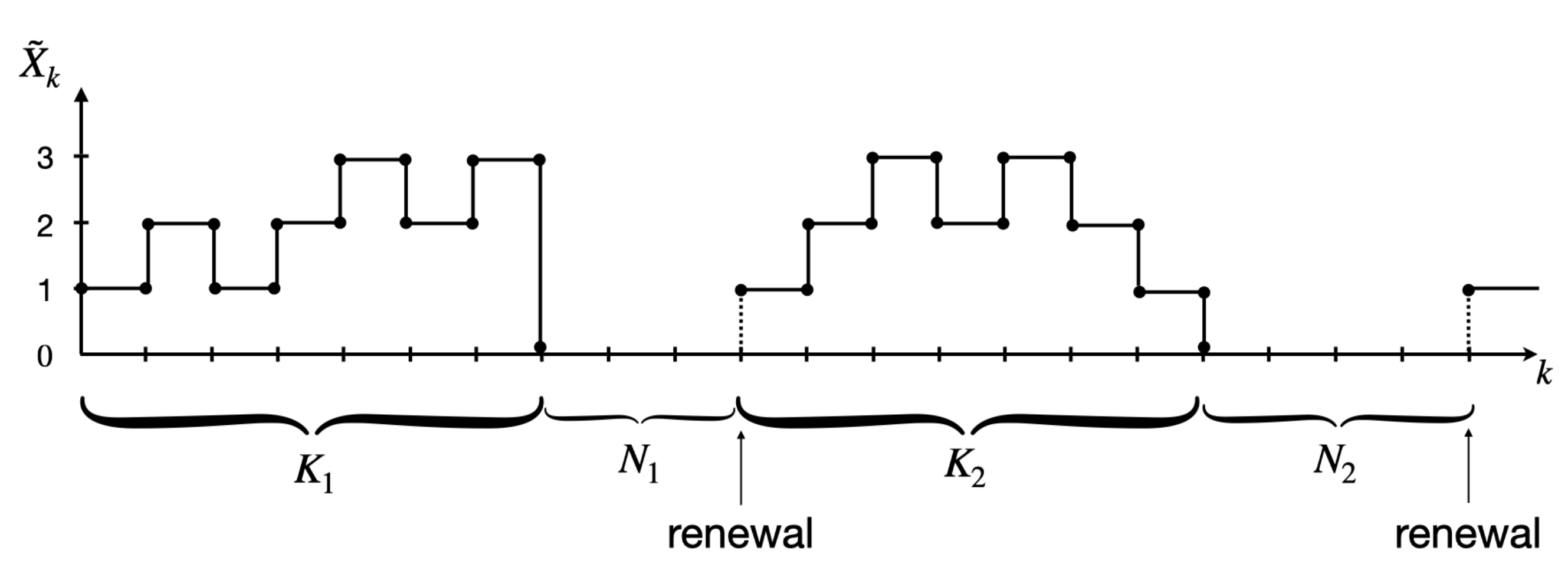}
        \caption{Illustration of the extended process $\{\Tilde{X}_k\}$. At first, the original process $X$ hits 0 at the $K_1$th jump (in total $X^{+}=4$ up-jumps); then stays at zero for $N_1$ number of time-steps until jumping up to 1, so a renewal occurs. Next, the process restarts and makes $K_2$ number of jumps until absorbed in 0, and then after $N_2$ time steps, a second renewal happens. $K_{1}$, $K_{2}$ are i.i.d. copies of $K$, and $N_{1},$ $N_{2}$ i.i.d. copies of $N.$}
    \label{fig:renewal_process}
    \end{center}
\end{figure}

{Next, we compute $\mu^{(A)}_{c}$. We recall that the size of an app-using component jumps up by 1 if there is a new app-user-infection, down by 1 if an app-user recovers naturally, and goes to 0 if an app-user is diagnosed. When only considering the number of jumps the process makes, we can forget the jump rates and simply let each time step correspond to a jump. The true process starts at $X_0=1$ and jumps up and down, or down to 0 with probabilities $\beta\pi/(\beta\pi+\gamma+\delta)$, $\gamma/(\beta\pi+\gamma+\delta)$ and $\delta/(\beta\pi+\gamma+\delta)$ respectively, and the process is absorbed when hitting 0 for the first time $K$. Let $X^{+}$ be the number of up-jumps before being absorbed in 0, and then we have 
\begin{equation}\label{eq:mu_c_1}
    \mu^{(A)}_{c} = 1+\mathbb{E}[{X}^{+}].
\end{equation}}
In order to get an expression for $\mathbb{E}[{X}^{+}]$, we extend the process $X=\{X_k\}$ to $\Tilde{X}=\{\Tilde{X}_k\}$ which is a renewal process (see Figure \ref{fig:renewal_process} for an illustration). More precisely, $ \Tilde{X}_k = X_k$ for $k \leq K$; once absorbed, we continue the $\Tilde{X}$ process by making it jump up to 1 with probability $\beta\pi/(\beta\pi+\gamma+\delta)$ and with the remaining probability $\Tilde{X}$ stays at 0. As a consequence, $\Tilde{X}$ will jump up to 1 after a $N\sim Geo(\beta\pi/(\beta\pi+\gamma+\delta))$ number of time steps ($N\ge 1$), and then the process restarts: a renewal has occurred. We also note that the process $\Tilde{X}$ makes an up-jump with probability $\beta\pi/(\beta\pi+\gamma+\delta)$ each time, whether currently at 0 or not. 

Let $\Tilde{X}^{+}_k$ denote the total number of up-jumps the extended process $\Tilde{X}$ has done up to $k'$th jump. Using the \textit{renewal reward theorem} (see e.g. Chapter 7.4 in \cite{ross2014introduction}), we have 
\begin{equation}\label{eq:renewal}
    \lim_{k \to \infty}\frac{\mathbb{E}[\Tilde{X}^{+}_k]}{k} = \frac{\mathbb{E}[{X}^{+}]+1}{\mathbb{E}[K]+\mathbb{E}[N]}.
\end{equation}
Clearly $\Tilde{X}^{+}_k\sim Bin(k,\ \beta\pi/(\beta\pi+\gamma+\delta))$ with expectation $k\beta\pi/(\beta\pi+\gamma+\delta).$ The left hand side of (\ref{eq:renewal}) hence equals $\beta\pi/(\beta\pi+\gamma+\delta)$. And as stated above, $N$ is geometrically distributed with success probability ${\beta\pi}/({\beta\pi+\gamma+\delta})$ on the support $\{1,2,3,...\}$, so
\begin{equation*}
    \mathbb{E}[N] = \frac{\beta\pi+ \gamma+\delta}{\beta\pi}.
\end{equation*}
It follows from (\ref{eq:renewal}) that
\begin{equation}\label{eq:mu_c_2}
    \mathbb{E}[{X}^{+}]+1 =1+\frac{\beta \pi}{\beta \pi+ \gamma+\delta} \mathbb{E}[K].
\end{equation}
To derive $\mathbb{E}[K]$, we first recall that $f(\pi)$ is the total number of non-app-users infected by an app-using component, so we can write $f(\pi)$ as  
\begin{equation*}
    f(\pi) = \mathbb{E}\left[\sum_{i=1}^{K} X^{(na)}_{i}\right],
\end{equation*}
where $X^{(na)}_{i}$ denotes the number of non-app-users infected between the $(i-1)$th and $i$th jump. For the distribution of $X^{(na)}_{i}$, we briefly state the idea as follows. At $(i-1)$th jump, suppose that there are $l$ infectives in the component, then the time until the next $i$th jump is exponentially distributed with intensity $l(\beta\pi+\gamma+\delta).$ During this time, such a component infects non-app-users at rate $l\beta(1-\pi)$. It follows that the distribution of $X^{(na)}_{i}$ is geometrically distributed with success probability  $$\frac{l({\beta\pi+\gamma+\delta})}{l({\beta\pi+\gamma+\delta}+\beta(1-\pi))}=\frac{{\beta\pi+\gamma+\delta}}{{\beta+\gamma+\delta}}$$ 
on the support $\{0,1,2,...\}.$ Because the success probability is independent of the current state $l$ implies that $X^{(na)}_{i}$ are not only independent but also identically distributed. Consequently, we obtain 
\begin{equation}\label{eq:N_c}
    f(\pi) =\mathbb{E}[K]\frac{\beta(1-\pi)}{\beta\pi+\gamma+\delta}. 
\end{equation}
Together with (\ref{eq:mu_c_1}) and (\ref{eq:mu_c_2}) gives 
\begin{equation}\label{eq:mu_c}
    \mu^{(A)}_{c} =1+\frac{\pi}{1-\pi}f(\pi). 
\end{equation}
Finally, it follows with (\ref{eq:expression_R_App_1}) that 
\begin{equation}\label{eq:expression_R_App}
    R^{(ind)}_{D-A} = \frac{f(\pi)}{1-\pi +\pi f(\pi)}.
\end{equation}

For non-app-users, it is easier to analyze since a given infectious non-app-user, on average, infects $(m_{21}+m_{22})$ individuals (app-users and non-app-users). Using (\ref{eq:expression_m21}) and (\ref{eq:expression_m22}) it follows that
\begin{equation}\label{eq:expression_R_NApp}
    R^{(ind)}_{D-N} = m_{21} + m_{22}=\frac{\beta}{\delta+\gamma},
\end{equation}
which is identical to $R_{0}$ in (\ref{eq:expression_R_0}). This is not surprising since non-app-users will never be contact traced. Plugging in (\ref{eq:expression_R_App}) and (\ref{eq:expression_R_NApp}) into (\ref{eq:expression_R_ind_pi}) gives 
\begin{equation}\label{eq:expression_R_ind_D}
 R^{(ind)}_{D}= \frac{\pi f(\pi)}{1-\pi +\pi f(\pi)} + (1-\pi) R_{0}.
\end{equation}
with $f(\pi)$ defined in (\ref{eq:expression_f}) and $R_{0}={\beta}/({\delta+\gamma})$. 


\subsubsection{Proof of Theorem \ref{thm:R_D_ind_monotone}.}
    Our aim is to prove that the derivative $R^{(ind)}_{D}{'}(\pi)$ is negative. Without loss of generality, we assume that $\delta+\gamma=1$, then we have $R_0 = \beta$ and 
    $f(\pi)=(1-\pi)g(\pi)/{2 \pi \delta}$
    with $g(\pi):=\beta \pi-1+\sqrt{(\beta \pi-1)^2+4\beta \pi \delta} \geq 0.$
   It implies with (\ref{eq:expression_R_ind_D}) that 
     \begin{align*}
        R^{(ind)}_{D}(\pi)&= \frac{(1-\pi)g(\pi)/2 \delta}{(1-\pi) + (1-\pi)g(\pi)/2 \delta} + (1-\pi) \beta 
        =  \frac{g(\pi)}{2\delta+g(\pi)}+ (1-\pi) \beta.
    \end{align*}
By differentiating the equation above w.r.t. $\pi$, we get 
\begin{align*}
    R^{(ind)}_{D}{'}(\pi) &= \frac{g{'}(\pi)(2\delta+g(\pi))-g(\pi)g{'}(\pi)}{(2\delta+g(\pi))^2}-\beta = \frac{2\delta g{'}(\pi)}{(2\delta+g(\pi))^2}-\beta.
\end{align*}
Further, the derivative of $g$ is given by
\begin{align*}
    g{'}(\pi)&=\beta + \frac{(\beta \pi -1)\beta+2\beta\delta}{\sqrt{(\beta \pi-1)^2+4\beta \pi \delta}}\\
    & = \beta  (1+\frac{\beta \pi-1+2\delta}{\sqrt{(\beta \pi-1)^2+4\beta \pi \delta}})\\
    & = \frac{\beta (2\delta+g(\pi))}{\sqrt{(\beta \pi+1)^2-4\beta \pi \gamma}}.
\end{align*}
Then $R^{(ind)}_{D}{'}(\pi)$ satisfies
\begin{align*}
    R^{(ind)}_{D}{'}(\pi) &= \frac{2\delta \beta (2\delta+g(\pi))}{(2\delta+g(\pi))^2\sqrt{(\beta \pi-1)^2+4\beta \pi \delta}} - \beta\\
    &= \frac{2\delta \beta}{(2\delta+g(\pi))\sqrt{(\beta \pi-1)^2+4\beta \pi \delta}}-\beta\\
    & = \frac{\beta h(\pi)}{(2\delta+g(\pi))\sqrt{(\beta \pi-1)^2+4\beta \pi \delta}}
\end{align*}
with 
\begin{align*}
   h(\pi)&:=2\delta-2\delta\sqrt{(\beta \pi-1)^2+4\beta \pi \delta}-g(\pi)\sqrt{(\beta \pi-1)^2+4\beta \pi \delta}\\
   &= 2\delta-2\delta\sqrt{(\beta \pi-1)^2+4\beta \pi \delta} - (\beta\pi -1)\sqrt{(\beta \pi-1)^2+4\beta \pi \delta}-(\beta \pi-1)^2-4\beta \pi \delta.
\end{align*}
Since $\beta, \delta >0$ and $g(\pi) \geq 0,$ we complete our proof by showing that 
\begin{align*}
    h(\pi) 
    &\leq 2\delta-2\delta(1-\beta\pi)-(\beta \pi-1)^2-(\beta \pi-1)^2-4\beta \pi \delta \\
    &= -2(\beta \pi-1)^2-2\beta\pi\delta <0,
\end{align*}
where the first inequality follows from that $\sqrt{(\beta \pi-1)^2+4\beta \pi \delta}\geq 1-\beta\pi.$

\subsubsection{Proof of Theorem \ref{thm:R_D_ind}.}
First, we assume $0<\pi<1.$ By multiplying both sides in (\ref{eq:expression_R_ind_D}) by $(1-\pi +\pi f(\pi))$,  $R^{(ind)}_{D} = 1$ is equivalent to 
\begin{equation*}
    \pi f(\pi)R_{0}(1-\pi)= (1-\pi)-(1-\pi)^2 R_{0}.
\end{equation*}
As a consequence, 
\begin{equation}\label{eq:result_f_pi}
    f(\pi)= \frac{1-R_{0}(1-\pi)}{\pi R_{0}}. 
\end{equation}
On the other hand, it follows from (\ref{eq:R_D}) that $R_{D} = 1$ is equivalent to 
\begin{equation*}
  \frac{( R_{0}(1-\pi))^2}{4}-R_{0}(1-\pi)+1 =\frac{( R_{0}(1-\pi))^2}{4}+f(\pi)R_{0}\pi
\end{equation*}
which gives the same result as in \eqref{eq:result_f_pi}.
Hence, we have showed that $R^{(ind)}_{D} =1$ is equivalent to $R_{D}=1$. Identical arguments can be used to show that $R^{(ind)}_{D} <1$ is equivalent to $R_{D} < 1$ as well as the opposite inequality.

In addition, when $\pi=0,$ it can be easily shown that $R^{(ind)}_{D}=R_{D}=R_{0}.$ If $\pi=1,$ i.e. everyone is an app-user, there is one app-using component that will eventually die out, implying that $R^{(ind)}_{D}$ will be below 1 (sub-critical). And using (\ref{eq:R_D}) gives $R_{D}=0.$ So both $R_{D}$ and $R^{(ind)}_{D}$ are below 1 when $\pi=1.$

\begin{rem}\label{remark:R_ind_D}
Similar to the situation for the individual component reproduction for manual contact tracing  (see \cite{zhang_analysing_2021}), this individual reproduction number $R^{(ind)}_{D}$ does not possess the traditional interpretation of the basic reproduction number: the average number of infections caused by infected people at the beginning of the outbreak. This is due to the delicate timing of event issues in an exponentially growing situation, closely related to those explained in \cite{ball_trapman_2016}. More specifically, when we derive the component reproduction number (upon which the individual reproduction number is based), we consider how many new components a component infects, but we do not keep track of when in time individuals of components get infected. For instance, root-individuals of components might on average infect more new individuals than other individuals in the component since they are not traced by their infector, \emph{and} root individuals are infected earlier than other individuals in the component. So, in the \emph{beginning} of an epidemic that grows exponentially, these root-individuals are over-represented when computing the mean number of infections among those infected early in the epidemic. However, our individual reproduction number $R^{(ind)}_{D}$ still possesses the important threshold property as {shown in Theorem \ref{thm:R_D_ind}}.
\end{rem}

\section{{Approximation of the initial phase in a large community for the combined model}}\label{sec:early_comb}
We now analyze the early stage of the epidemic with both manual and digital contact tracing. First, we define the limit process $E_{DM}(\beta,\gamma,\delta,p,\pi)$ as follows. In the process, individuals are either app-users or non-app-users. There is one initial ancestor who has the same type of the initial case in $E^{(n)}_{DM}(\beta,\gamma,\delta,p,\pi)$. During the lifetime, each individual gives birth to app-users with manual tracing link at rate $\beta\pi p$, and without link at rate $\beta\pi (1-p)$, and to non-app-users with manual tracing link at rate $\beta(1-\pi)p$, and without link at rate $\beta(1-\pi)(1-p).$ Individuals recover naturally at rate $\gamma$ and are diagnosed at rate $\delta.$ Once diagnosed, each of its descendants and parent with manual tracing link are diagnosed immediately, and so on. If such a diagnosed individual is an app-user, all of its app-using descendants as well as its parent (in case of an app-user), will also be diagnosed at the same time, and so on. Please note that tracing is also iterated for traced individuals who have recovered naturally upon diagnosis, so their parents and descendants are diagnosed according to the same rule. We first prove Proposition \ref{prop:earlyapprox_comb}, which is similar to the proof of Proposition \ref{prop:earlyapprox_dct}.

\subsection{{Proof of Proposition \ref{prop:earlyapprox_comb}}}\label{sec:proof_thm_comb}
First, we note that the only difference between the $n$-epidemic $E^{(n)}_{DM}(\beta,\gamma,\delta,p,\pi)$ and the limit process $E_{DM}(\beta,\gamma,\delta,p,\pi)$ is that the rate of new infections in the epidemic is deflated by the factor $S(t)/n$ since infectious contacts only result in infection if the contacted person is still susceptible. Using standard coupling arguments as described in the proof of  Proposition \ref{prop:earlyapprox_dct}, it can be shown that the time $T_n$ of the first contact with an already infected person tends to $\infty$ in probability as $n\to\infty$, and prior to this the epidemic and the limiting process agree completely. Consequently, the $n$-epidemic converges in distribution to $E_{DM}(\beta,\gamma,\delta,p,\pi)$ on any finite time interval.
     
\begin{figure}[tb!]
    \begin{center}
        \includegraphics[width=80mm]{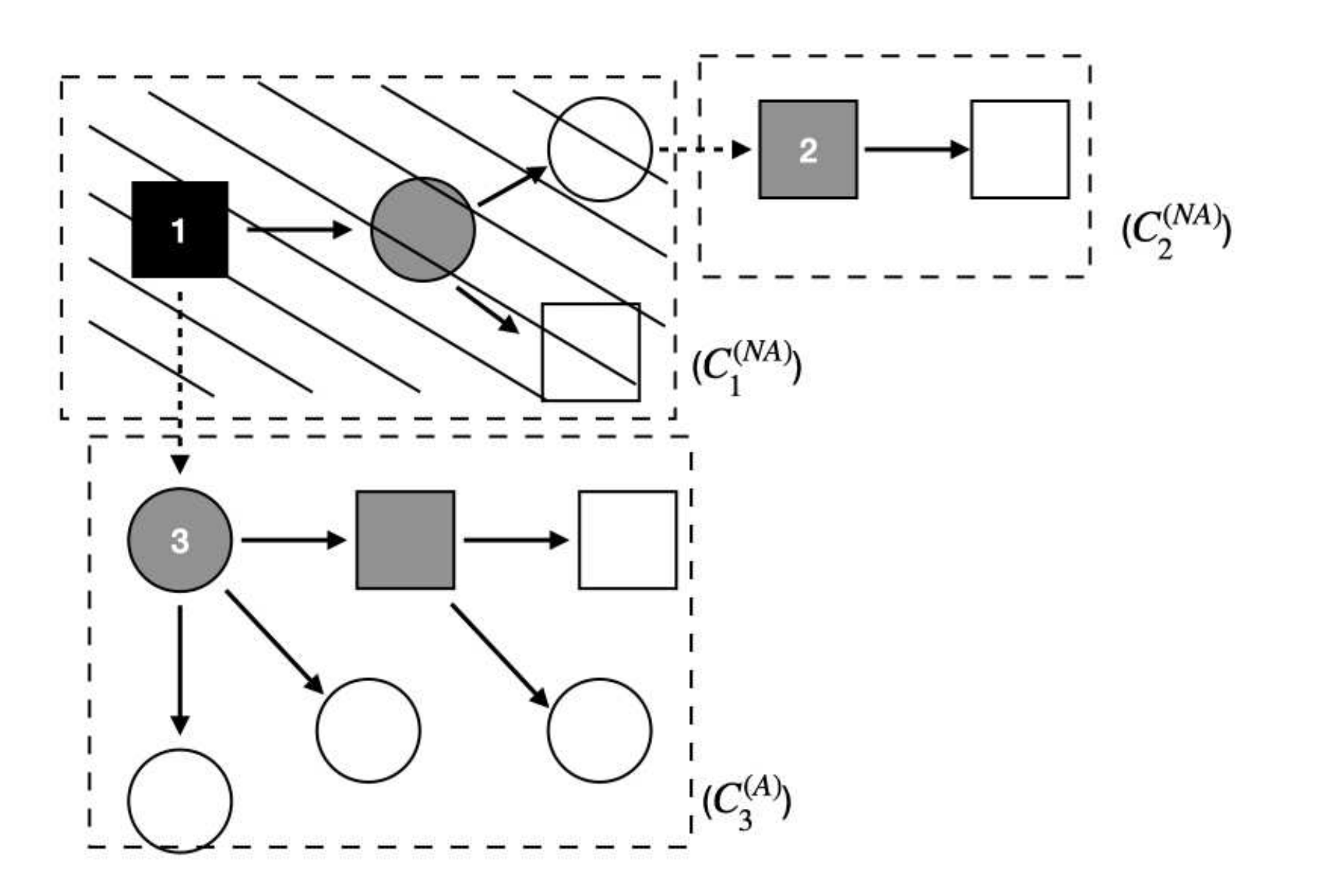}
    \caption{Example of the process of to-be-traced components started by app-users and non-app-users for the combined model: the circles are app-users and squares non-app-users. The nodes in white, black, and gray are ``infectious", ``diagnosed" and ``naturally recovered", respectively. The rectangular region surrounded by the dashed line is for the component, whereas the whole rectangle is dashed when the component is diagnosed. Solid edges stand for contacts that could be traced either by manual or digital contact tracing (so always full edges between each pair of app-users), whereas dashed edges are infectious contacts that will not be traced.}
    \label{fig:two_type_br_combine}
    \end{center}
\end{figure}

Following a similar idea as in the proof of Theorem \ref{sec:proof_thm_dct} in Section \ref{sec:proof_thm_dct}, we describe this limit process $E_{DM}(\beta,\gamma,\delta,p,\pi)$ in terms of to-be-traced components. For the combined model, the app-users could not only be traced by app-users, but also by non-app-users with probability $p.$ Further, non-app-users could be traced by individuals (either app-users or non-app-users) with probability $p$. Accordingly, to-be-traced components also contain infections between app-users and non-app-users and pairs of app-users, who will be traced manually. Moreover, new components are created if an app-user infects a non-app-user without manual contact tracing taking place or if a non-app user infects either type and there is no manual contact tracing. The new components hence differ only in how they are created: with an app-user (infected by a non-app-user without manual contact tracing) or a non-app-user (without manual contact tracing), and given the current state, the evolvement of such components does not depend on how they were initiated and is hence the same for both types of components. There will be no contact tracing between any of these components, which thus makes them evolve independently, as stipulated for branching processes. Here, we give an example in Figure \ref{fig:two_type_br_combine} to illustrate how the two types of components (starting with an app-user or a non-app-user) grow and reproduce. We start with case 1, which is the root of the component $C^{(NA)}_{1}$, which generates two new components $C^{(NA)}_{2}$ and $C^{(A)}_{3}$. The component $C^{(NA)}_{2}$ is  with root 2 as a non-app-user and $C^{(A)}_{3}$ starts with an app-user 3. Once individual 1 is diagnosed, the whole component $C^{(NA)}_{1}$ is diagnosed.

\subsection{Proof of Corollary~\ref{col:combined}}
We see that the combined model can be approximated by a suitable two-type branching process during the early stage of the epidemic, where type-1 ``individuals" are to-be-traced components starting with an app-user and type-2 are with a non-app-user as root. Let $M^{(c)}$ be the corresponding mean offspring matrix
\begin{equation}\label{eq:M_c}
    M^{(c)}=\begin{pmatrix}
m^{(c)}_{11} & m^{(c)}_{12}\\\\
m^{(c)}_{21}  & m^{(c)}_{22}
\end{pmatrix},
\end{equation} where element $m^{(c)}_{ij}$, is the mean number of new components of type $j$ produced by a single component of type $i$, for $i,j =1,2.$ 
Let $R_{DM}$ be the largest eigenvalue of $M^{(c)}$, then $R_{DM}$ can be expressed similarly to the digital reproduction number $R_{D}$ as shown in (\ref{eq:eigenvalue_M}), but with $m^{(c)}_{ij}$ replacing $m_{ij}$.
Similar to $R_{D}$,  $R_{DM}$ determines whether a major outbreak is possible or not. We refer $R_{DM}$ to the \textit{effective component reproduction number for the combined model}.
As shown in the proof of Corollary~\ref{col:digital}, it follows by the standard multi-type branching process theory \cite{ref2,becker_multi_type_1990} that the limit process $E_{DM}$ will surely go extinct if $R_{DM} \leq 1$. If $R_{DM}>1$, the process will grow beyond all limits with non-zero probability. Then, by Proposition \ref{prop:earlyapprox_comb}, it can be shown that if $R_{DM} > 1$, there will be a large outbreak of the epidemic with non-zero probability, whereas if $R_{DM} \leq 1,$  there will be a minor outbreak for sure.

In the following text, we derive the elements $m^{(c)}_{ij}$. We first note that the to-be-traced components may also contain non-app-users that are manually traced. So, it is necessary to track both the number of infectious app-users and infectious non-app-users in the component. Let $N^{(i)}_{j}(t)$ be the number of infectious type-$j$ individuals in the component starting with one type-$i$ individual at time $t$ $(i,j=1,2),$ where we set type-$1$ to be app-users and type-$2$ to be non-app-users. When the components are born ($t=0$), we have $N^{(1)}(0)=(1,0)$, whereas  $N^{(2)}(0)=(0,1)$. Suppose that at time $t$ there are $k$ infectious app-users and $l$ infectious non-app-users in a component starting with a type-$i$ individual, denoted by
\begin{equation*}
    N^{(i)}(t):= (N^{(i)}_{1}(t), N^{(i)}_{2}(t))=(k,l).
\end{equation*}
There are 5 different competing events that can occur: an app-user is newly infected by an app-user or with a manual tracing link by a non-app-user ($N^{(i)}(t)+(1,0)$), one of the $k$ app-users recovers naturally ($N^{(i)}(t)-(1,0)$), a non-app-user is newly infected with a manual tracing link irrespective of the type of the infector ($N^{(i)}(t)+(0,1)$), one of the $l$ non-app-users recovers naturally ($N^{(i)}(t)-(0,1)$), or one of the $(k+l)$ infectious individuals (of either type) is diagnosed, so the whole component is diagnosed (drops to $(0,0)$). The corresponding transition rates are given by
$$
N^{(i)}(t^{-}) \to \left\{
\begin{array}{lll}
 N^{(i)}(t)+(1,0)    &  \text{at rate\ } k\beta\pi+l\beta\pi p \\
  N^{(i)}(t)+(-1,0)   & \text{at rate\ } k\gamma\\
  N^{(i)}(t)+(0,1) & \text{at rate\ } (k+l)\beta(1-\pi)p\\
  N^{(i)}(t)+(0,-1) & \text{at rate\ } l\gamma\\
  (0,0) & \text{at rate\ } (k+l)\delta .
\end{array}\right.
$$

Next, we describe the birth of new components. We note that only infections with no digital or manual contact tracing links create a new component. This can only occur to pairs of non-app-users or pairs of app-users and non-app-users, but it never happens between two app-users. As a consequence, given that $N^{(i)}(t)=(k,l)$ ($k$ infectious app-users and $l$ infectious non-app-users), the component gives birth to new components with an app-user as root at rate $l\beta\pi(1-p)$, and new components with a non-app-user as root at rate $(k+l)\beta(1-\pi)(1-p).$ Based on the previous discussion, we have derived the birth rates of new components of type $i$ conditional on $N^{(i)}(t)$. The expected total numbers  $m^{(c)}_{ij}$ of new type-$j$ components produced by a type-$i$ component are hence the expected values of these expressions integrated over time:
\begin{align}
\label{eq:elements_m_c_11}
    m^{(c)}_{11} &= \int_{0}^{\infty}\mathbb{E}[N^{(1)}_{2}(t)]\beta\pi(1-p) dt,\\
    \label{eq:elements_m_c_12}
    m^{(c)}_{12} &=  \int_{0}^{\infty}(\mathbb{E}[N^{(1)}_{1}(t)]+\mathbb{E}[N^{(1)}_{2}(t)])\beta(1-\pi)(1-p) dt,\\
    \label{eq:elements_m_c_21}
    m^{(c)}_{21} &= \int_{0}^{\infty}\mathbb{E}[N^{(2)}_{2}(t)]\beta\pi(1-p) dt,\\
    \label{eq:elements_m_c_22}
    m^{(c)}_{22} &= \int_{0}^{\infty}(\mathbb{E}[N^{(2)}_{1}(t)]+\mathbb{E}[N^{(2)}_{2}(t)])\beta(1-\pi)(1-p) dt.
\end{align}

We have not been able to derive $\mathbb{E}[N^{(i)}_{j}(t)]$ analytically. The main reason is that not all the four jump rates are linear in 
the number of infectious individuals. In particular, the rate at which $N^{(i)}=(k,l)$ increases in the first component equals $k\beta\pi+l\beta\pi p$ (non-linear in $k+l$).

 {In Section \ref{sec:numerical}, we perform simulations to derive the
 $m^{(c)}_{ij}$ and numerically compute the component reproduction number $R_{DM}$ for given parameters. Based on the numerical results, we analyze how $R_{DM}$ depends on the manual tracing probability $p$ and the app-using fraction $\pi$.}

\section{Numerical Investigation}\label{sec:numerical}

In this section, we numerically illustrate the effects of digital contact tracing, as well as the combined effect of both types of contact tracing. 
For all the numerical results below (except Table \ref{tab:compare_ct}), we set the epidemic parameters to be fixed: 
$\beta = \frac{6}{7}, \gamma =\frac{1}{7}$ so that the reproduction number without any testing or contact tracing is ${\beta}/{\gamma}=6$. When also considering diagnosis ($\delta >0$), we would have the basic reproduction number $R_0 = {\beta}/({\gamma+\delta})$ smaller than 6. The code used for this section is available on GitHub at \cite{githublink}.

\subsection{Comparison of digital and manual tracing only}

\begin{figure}[tb!]
\begin{subfigure}[tb!]{0.48\textwidth}
    \begin{center}
        \includegraphics[width=\textwidth]{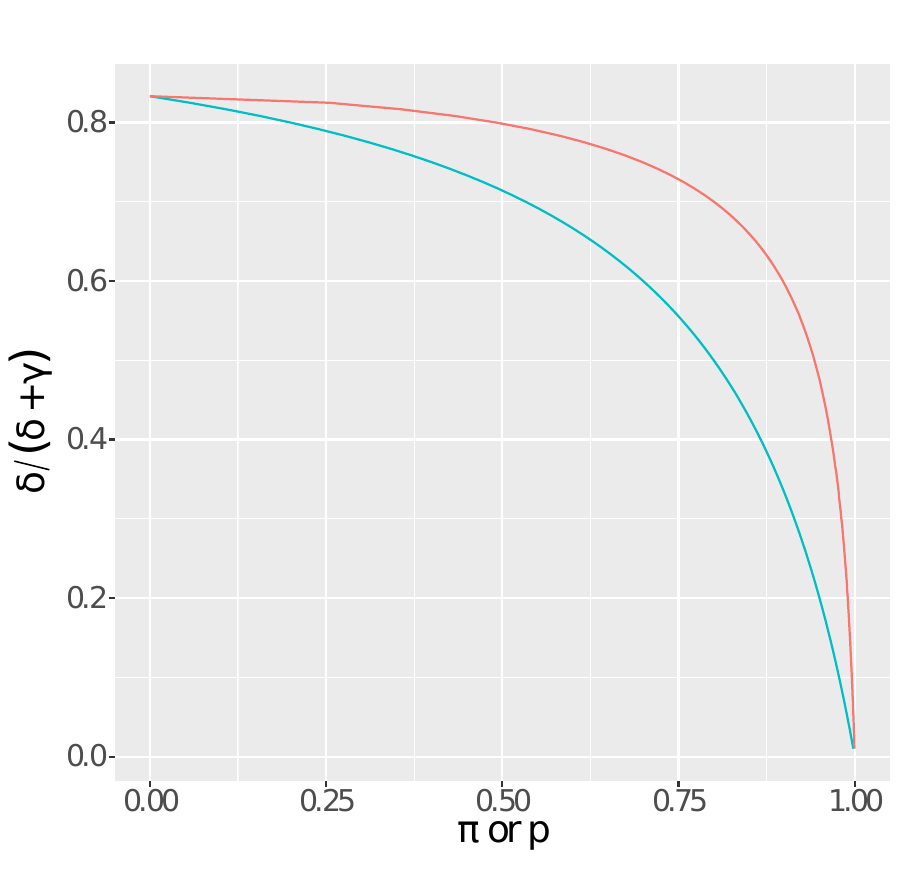}
    \caption{}
    \label{fig:compare_DM}
    \end{center}
    \end{subfigure}
\hfill
\begin{subfigure}[tb!]{0.56\textwidth}
\begin{center}
\includegraphics[width=\textwidth]{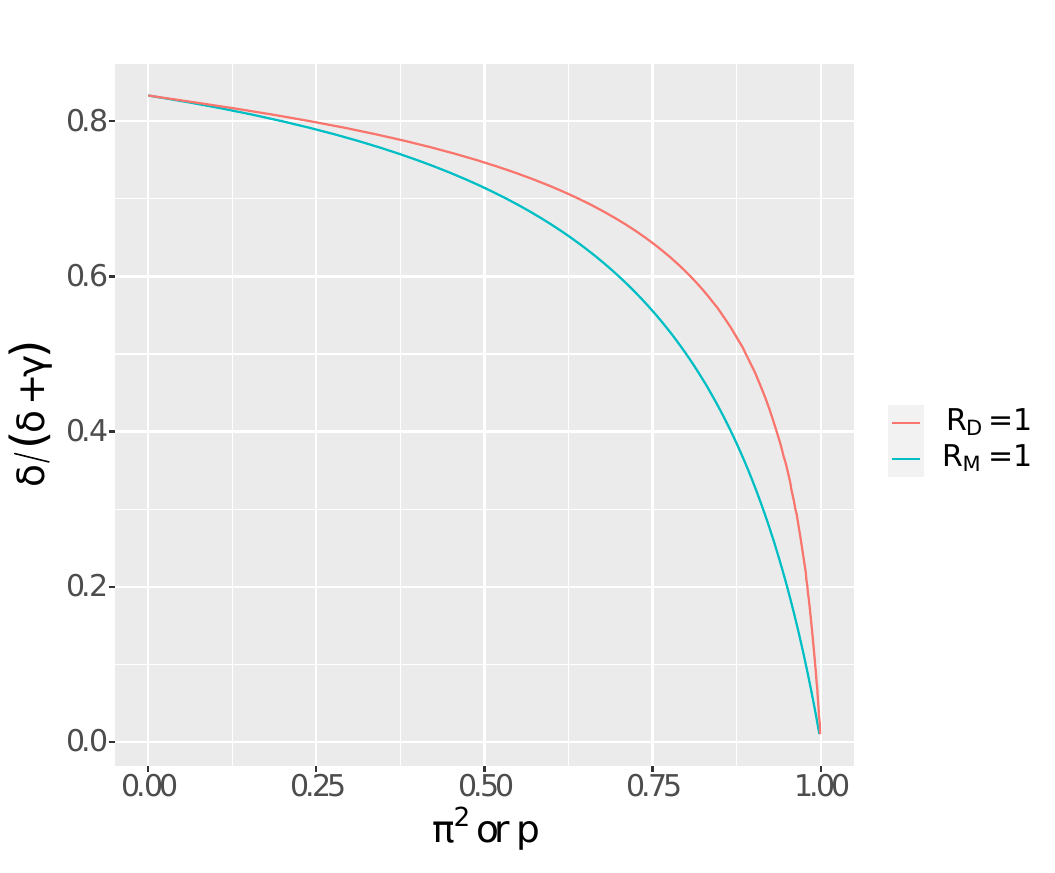}
\caption{}
\label{fig:compare_DM_sq}
\end{center}
\end{subfigure}
\caption{Plot a). The red curve shows which combinations of the testing fraction ${\delta}/({\delta + \gamma})$ and fraction $\pi$ of app-users for which $R_D=1$ for digital contact tracing, and the blue curve which combinations of testing fraction ${\delta}/({\delta + \gamma})$ and tracing probability $p$ for which $R_M=1$ in the manual contact tracing. Plot b) shows the same curves, but now having $\pi^2$ on the $x$-axis (motivated by both individuals having to be app-users for digital tracing to take place).}
\end{figure}

 In Figure \ref{fig:compare_DM}, we plot the necessary testing rate measured by testing level ${\delta}/({\delta + \gamma})$ for which the manual $R^{(c)}_{M}$ given by (\ref{eq:expression_R_M}) equals 1  and digital component reproduction number $R^{}_{D}$ given by (\ref{eq:R_D}) equals 1, for a given value of $p$ and $\pi$ respectively. {The figure confirms Theorem \ref{thm:p_c_pi_c} that} that for a given testing level ${\delta}/({\delta + \gamma})$, a larger app-using fraction $\pi$ is required compared to the tracing probability $p$ for the model with manual contact tracing. 
 
 One potential explanation for this could be that the effectiveness of digital contact tracing is more related to the square of the fraction of app-users ($\pi^2$): the digital tracing is carried out only if both the confirmed case and the contact are app-users. In Figure \ref{fig:compare_DM_sq}, we therefore plot the critical combination of ${\delta}/({\delta + \gamma})$ against $\pi^2$ instead. It can be seen that the curve for $\pi^2$ is closer to the one for $p$ but still different: we still need a larger value on $\pi^2$ for the reproduction number to equal unity. An alternative explanation would be that there is a fraction $(1-\pi)$ of individuals will never be digitally contact traced, whereas manual tracing can reach all individuals. This suggests that digital tracing requires a higher effectiveness among app-users to compensate for the absence of contact tracing among non-app-users.

\subsection{The individual and component reproduction numbers for digital tracing only}

\begin{figure}[tb!]
\begin{subfigure}[tb!]{0.48\textwidth}
\begin{center}
    \includegraphics[width=\textwidth]{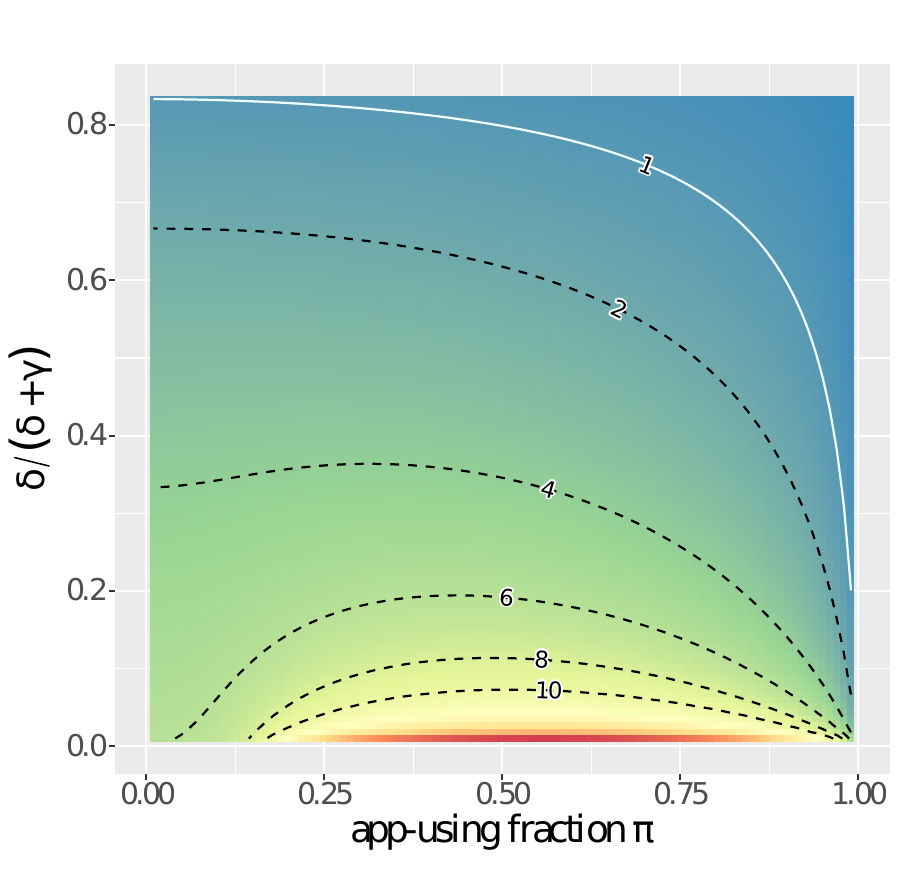}
\caption{Heatmap of $R_{D}$}
\label{fig:RD}
\end{center}
\end{subfigure}
     \hfill
\begin{subfigure}[tb!]{0.6\textwidth}
\begin{center}
\includegraphics[width=\textwidth]{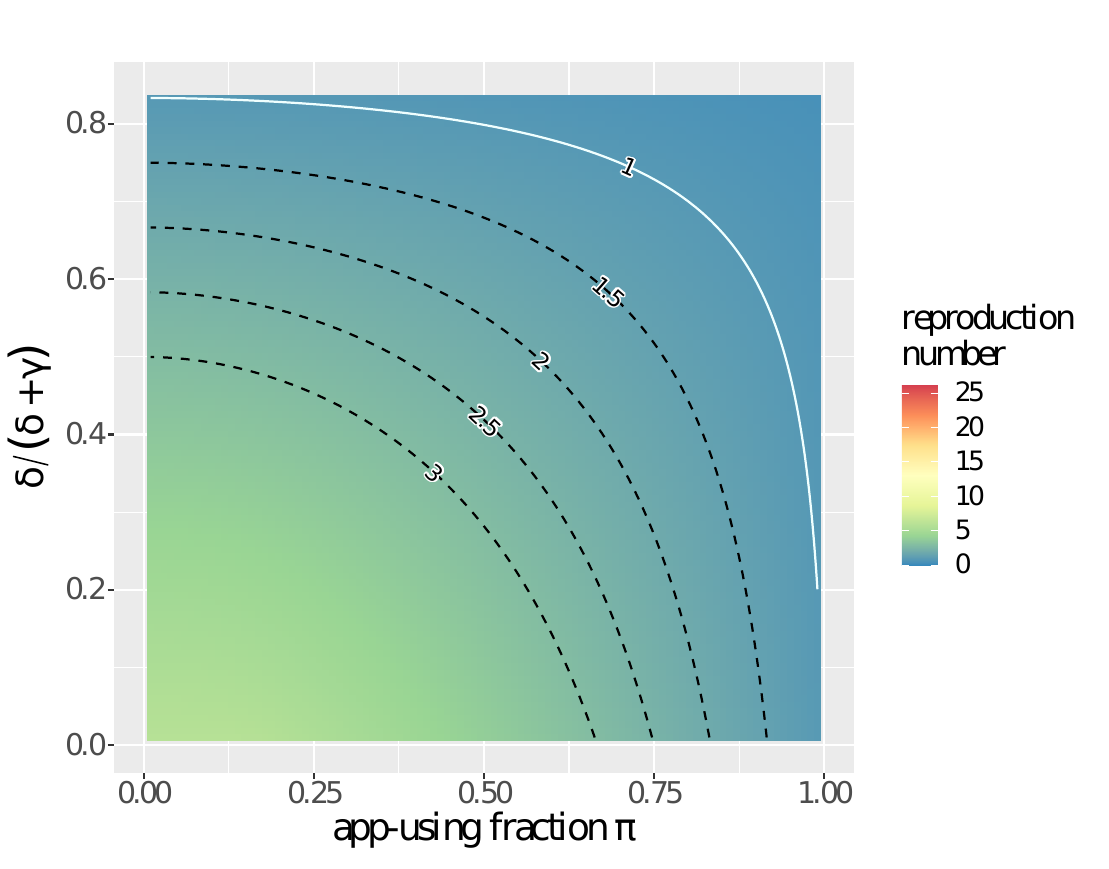}
\caption{Heatmap of $R^{(ind)}_{D}$}
\label{fig:RDind}
\end{center}
\end{subfigure}
\caption{Heatmaps of the effective reproduction numbers for digital tracing: \ref{fig:RD} for $R_{D}$ and \ref{fig:RDind} for $R^{(ind)}_{D}$ varies with testing fraction ${\delta}/({\delta + \gamma})$ in $[0.01,0.83]$ and fraction of app-users $\pi$ in $[0.01,0.99]$; $\beta=\frac{6}{7}, \gamma=\frac{1}{7}$ fixed.  The white solid lines indicate where the corresponding reproduction number equals 1, the black dashed lines in \ref{fig:RD} are for $R_{D}=2,4,6,8,10$ from up to down and in \ref{fig:RDind} for the cases when $R^{(ind)}_{D}$ equals $3, 2.5,2,1.5,$ respectively from left to right.}

\label{fig:R_ct}
\end{figure}

In Figure \ref{fig:RD} and \ref{fig:RDind}, we plot how the digital reproduction number $R_{D}$ and individual reproduction number $R^{(ind)}_{D}$ vary with the testing fraction ${\delta}/({\delta + \gamma})$  and the fraction of app-users $\pi$.  By looking at
the colors at the bottom {or the dashed curves at $R_{D}=6, 8,$ and 10} in Figure \ref{fig:RD}, we note that the component reproduction number $R_{D}$ is not monotonically decreasing in $\pi$. This is at first surprising since we expect a smaller reproduction the more app-users there are.

This non-monotonicity of $R_{D}$ is explained as follows. First of all, it can be easily seen from (\ref{eq:expression_m21}) and (\ref{eq:expression_m22}) that the element $m_{21}$ is increasing and $m_{22}$ is decreasing with $\pi.$ Then, at a very low level of testing fraction, increasing $\pi$ leads to more app-users, and thus the \emph{size} of the app-using component will grow even if fewer non-app-users will be infected \emph{per} app-user. The average number of non-app-users infected by one app-using component $m_{12}$ may hence increase with $\pi$. Then according to the expression for $R_{D}$ in (\ref{eq:R_D}), $R_{D}$ could increase with $\pi$ when the testing fraction ${\delta}/({\delta + \gamma})$ is low. On the other hand, we see from Figure \ref{fig:RDind} that the individual reproduction number $R^{(ind)}_{D}$ \textit{is} monotonically decreasing in $\pi$ as was proven in Theorem \ref{thm:R_D_ind_monotone}.

\begin{figure}[tb!]
    \begin{center}
        \includegraphics[width=.9\textwidth]{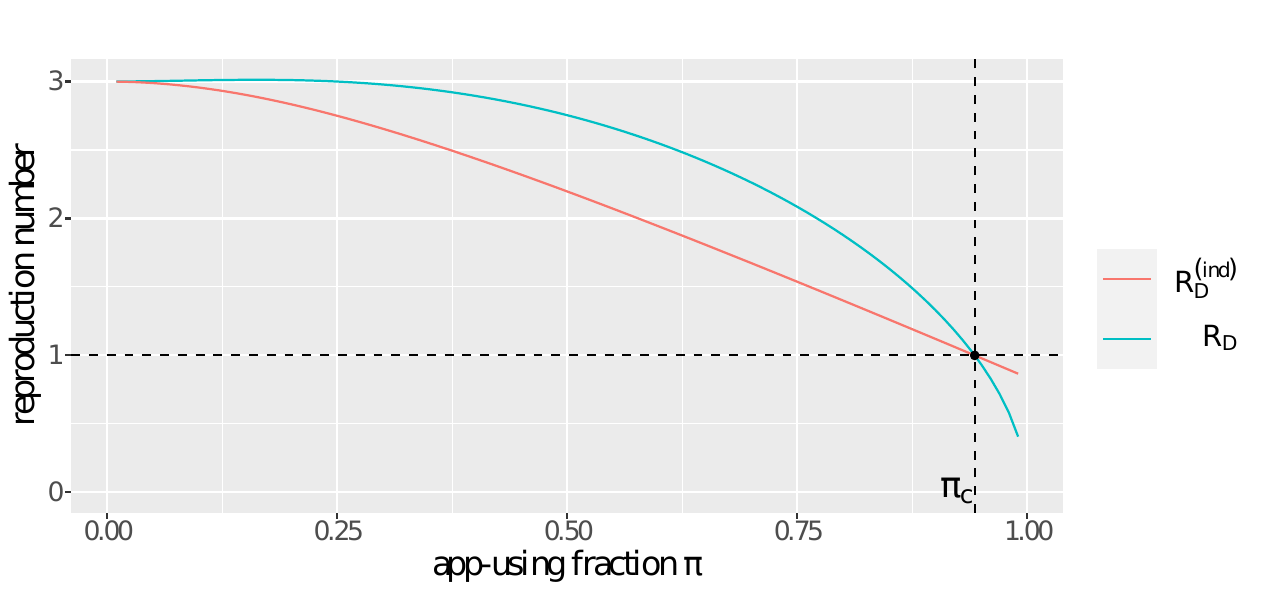}
    \caption{Plot of two curves of $R^{(ind)}_{D}$ (in red) and $R^{}_{D}$ (in blue) against $\pi$ with $\beta=\frac{6}{7}, \gamma=\frac{1}{7}$ and $\delta=\frac{1}{7}$ fixed. The horizontal dashed line stands for when the reproduction number equals 1{, and the vertical dashed line for the critical $\pi_{c} \approx 0.94$.}}
    \label{fig:compare_ind_D}
    \end{center}
\end{figure}

In Figure \ref{fig:compare_ind_D} we illustrate Theorem \ref{thm:R_D_ind} stating that the two reproduction numbers $R^{}_{D}$ and $R^{(ind)}_{D}$ have the same threshold at 1. As before, we fix $\beta=\frac{6}{7}$, $\gamma = \delta = \frac{1}{7}$ such that $R_0=3$ and half of the infected are diagnosed (and subject to being contact traced). The two reproduction numbers are plotted as a function $\pi$, the fraction of app-users. The figure shows that the reproduction numbers are different, but they both attain the value 1 for the same critical fraction of app-users $\pi_c \approx 94\%$.

\subsection{Effect of combining digital and manual tracing}

\begin{figure}[tb!]
\begin{subfigure}[tb!]{0.48\textwidth}
\begin{center}
    \includegraphics[width=\textwidth]{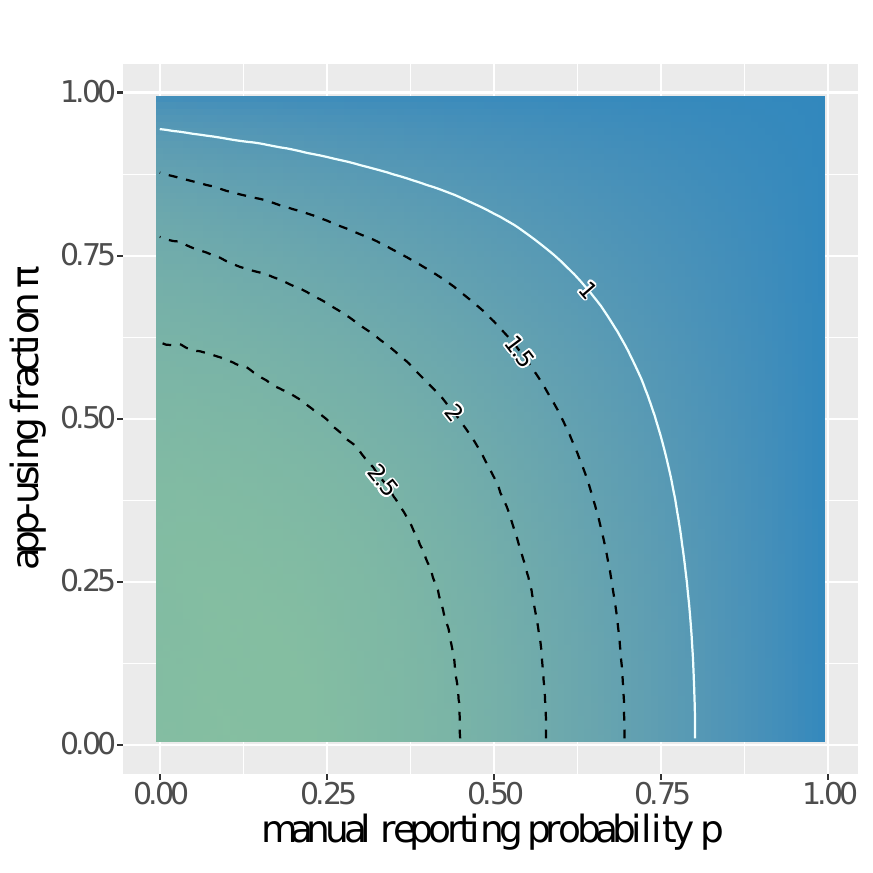}
\caption{Heatmap of $R_{DM}$ with $\frac{\delta}{\delta+\gamma}=0.5$}
\label{fig:RDM}
\end{center}

\end{subfigure}
     \hspace{4mm}
\begin{subfigure}[tb!]{0.6\textwidth}
\begin{center}
\includegraphics[width=\textwidth]{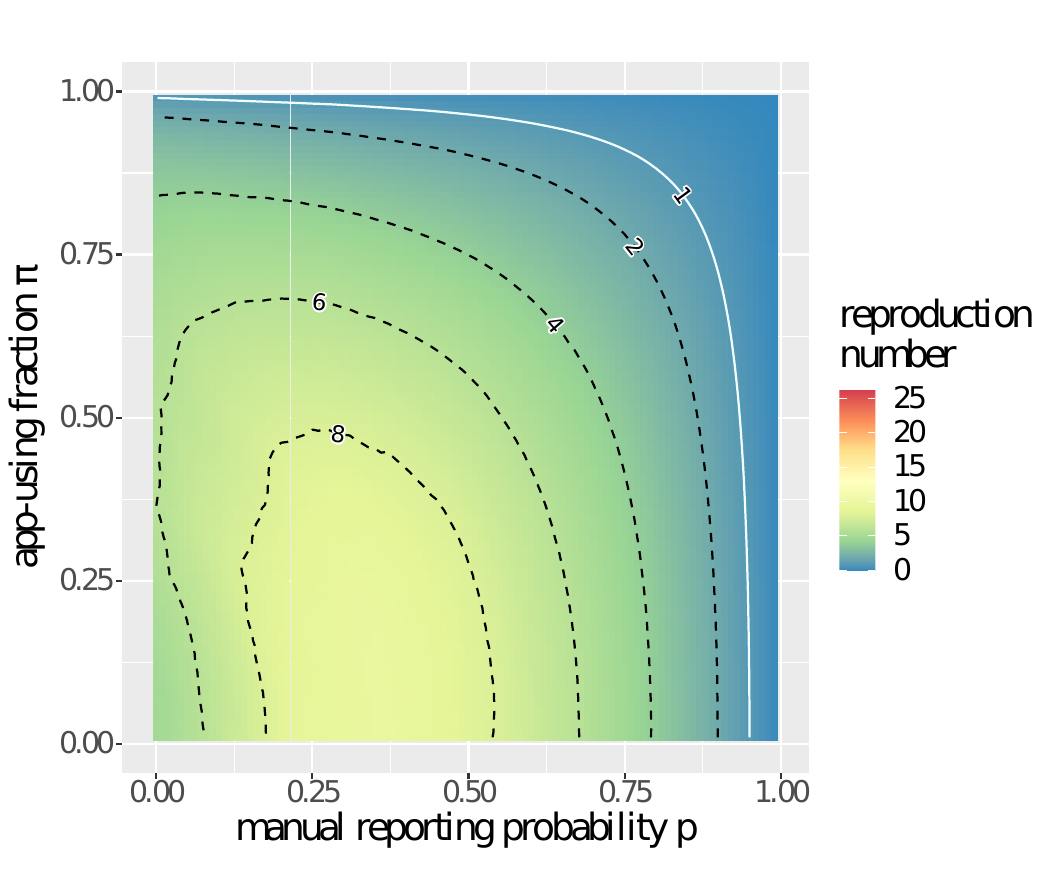}
\caption{Heatmap of $R_{DM}$ with $\frac{\delta}{\delta+\gamma}=0.2$}
\label{fig:RDMsmall}
\end{center}
\end{subfigure}
\caption{{Heatmaps of the effective reproduction number $R_{DM}$ varying with manual tracing probability $p$ in $[0,0.99)$ and fraction of app-users $\pi$ in $[0,0.99)$; with $\beta=\frac{6}{7}, \gamma=\frac{1}{7}$ fixed, \ref{fig:RDM} for $\delta=\frac{1}{7}$ and \ref{fig:RDMsmall} for smaller $\delta=\frac{1}{28}$. The white solid lines indicate where $R_{DM}=1$, the black dashed lines in \ref{fig:RDM} are for the cases when $R_{DM}$ equals $2.5,2,1.5,$ respectively from left to right, and in \ref{fig:RDMsmall} for $R_{DM}=8,6,4,2$.}}

\label{fig:R_dm}
\end{figure}

We now illustrate the effect of the combined model by computing $R_{DM}$ numerically as follows. For $\beta=\frac{6}{7}, \gamma=\frac{1}{7}$ and two distinct values for $\delta$: $\frac{1}{7}$ and $\frac{1}{28}$, and a grid of values for $\pi$ (the fraction of app-users) and $p$ (the fraction of contacts that are traced manually), we perform 10000 simulations to compute the integrals in (\ref{eq:elements_m_c_11})-(\ref{eq:elements_m_c_22}) for each simulation (without expectation). We then take averages over the simulations to obtain estimates of the four elements $m^{(c)}_{ij}$, and compute the largest eigenvalue to obtain estimates of $R^{}_{DM}$. In Figure \ref{fig:RDM} and \ref{fig:RDMsmall}, we show the heatmaps of $R^{}_{DM}$ as a function of  $p$ and $\pi$ for the two choices of $\delta$. We see from the left figure that $R^{}_{DM}$ seems monotonically decreasing in both $\pi$ and $p$ as expected when $\delta=\frac{1}{7}$ (as one would expect), but this fails when $\delta$ is smaller as seen in the right figure. This can be explained in a similar way as for $R_{D}$. For example, suppose $\pi$ and $\delta$ are both small, then the size of a component (either starting with an app-user or non-app-user) may increase with $p$ since more infections happen within a component. As a consequence, such a component may, in fact, infect slightly more new components, thus resulting in a larger reproduction number. However, on the individual scale, a much larger component infecting slightly more new components would still imply that each individual will infect fewer new individuals. Here, it might be possible to derive an individual reproduction number, $R^{(ind)}_{DM}$, which one would expect decreases monotonically in both $\pi$ and $p$, but we have not pursued this wrote, partly because the expression for the component reproduction number $R^{}_{DM}$ is complicated in itself. 
But if we focus on the curve $R_{DM}=1$ in Figure \ref{fig:RDMsmall}, this curve seems to be monotone, suggesting that decreasing one of the contact tracing procedures requires the other to be increased in order to maintain herd immunity. A proof of this is missing.
 
All results above are for the limiting processes $E_D$ and $E_{DM}$ respectively. To see how the results fit the finite $n$ case, we perform $10000$ simulations of the epidemic for $n = 1000$ and $5000$, starting with one initial infective. We choose the parameters in four different cases: without any contact tracing ($p=\pi=0$); with digital tracing only ($p=0,\pi=\frac{2}{3}$); with manual tracing only: ($p=\frac{2}{3},\pi=0$) and with both manual and digital tracing ($p=\pi=\frac{2}{3}$), while keeping $\beta=0.8, \gamma=\frac{1}{7}, \delta=\frac{1}{7}$ fixed. It can be seen from Table \ref{tab:compare_ct} that the fraction of major outbreaks (defined by having more than $10\%$ infected) is distinct away from 0 when the reproduction number exceeds 1 and is very close to 0 when the reproduction is smaller than 1 (thus supporting {Corollary \ref{col:digital} and \ref{col:combined}}). In addition to our epidemic simulations, we have included a column for the limiting results (denoted by $n=\infty$), where we simulate the corresponding limiting branching process with a single ancestor 10000 times. We then determine the limit fraction of major outbreaks defined as scenarios with more than 1000 individuals born. The results indicate that our limiting approximations are effective.

\begin{table}[tb!]
    \centering
    \caption{Results from 10000 simulated epidemics and the corresponding limiting branching processes}\label{tab:compare_ct}
    \footnotesize{
    \begin{tabular}{|l| l| l| l| l| l| l| l| l|}
    \hline
     &  & Reproduction & \multicolumn{3}{c|}{Fraction of} & \multicolumn{3}{c|}{Mean fraction of infected}\\
    $p$ & $\pi$  &    number   &   \multicolumn{3}{c|}{major outbreaks} &  \multicolumn{3}{c|}{among major outbreaks} \\
         \cline{4-9}
         &       &            & $n=1000$ & $n=5000$ & $n=\infty$  &   $n=1000$ & $n=5000$ &  \\
    \cline{1-8}
    0 & 0 & $R_{0}=2.80$ & 0.637 & 0.638 & 0.639 & 0.924 & 0.925 & \\
    0 & $\frac{2}{3}$ & $R_{D}=2.20$ & 0.478 & 0.482 & 0.480 &  0.812  &0.814 &$-$\\
    $\frac{2}{3}$ & 0 & $R_{M}=1.49$ & 0.280  & 0.276 & 0.282 & 0.492  & 0.507&\\
    $\frac{2}{3}$ & $\frac{2}{3}$ & $R_{DM}=0.92$ & 0.067 & 0.014 & 0.011 & 0.200  &0.147 &\\
    \hline
    \end{tabular}
    }
\end{table}

\section{Discussion}\label{sec:discussion}
In this paper, we analyzed a Markovian epidemic model with digital contact tracing, without or with manual contact tracing. For the epidemic model with digital tracing only, the early stage of the epidemic was proven to converge to a two-type branching process (as the population size $n\to\infty$)  with one type being to-be-traced app-using components and the other being non-app-users, and $R_D$ being the largest eigenvalue of the of mean offspring matrix of the two-type branching process. An individual reproduction number $R^{(ind)}_{D}$ was also derived. 
It was investigated both analytically and numerically how $R_D$ and $R^{(ind)}_{D}$ depended on model parameters as well as on each other. The reproduction number $R_D$ was shown not monotone in $\pi$ (the fraction of app-users) in general, whereas the individual reproduction number $R^{(ind)}_{D}$ was proven to be monotonically decreasing in $\pi$. When comparing digital tracing with manual tracing, it was proven that a smaller tracing probability $p$ is needed than the corresponding app-using fraction $\pi$ in order to lower the reproduction number down to 1. 

We then analyzed the combined model and proved that the beginning of the epidemic converges to a different two-type branching process, with both types being to-be-traced components and only differing on the type of the root (app-user or non-app-user). The corresponding reproduction number had a rather complicated expression, which currently can only be evaluated numerically, and it was shown numerically not to be monotonically decreasing in $\pi$ and $p$ in general. 

There are several extensions of the present model that would make it more realistic. On the one hand, we assume that there is no delay in both manual and digital contact tracing, i.e. any traced contact (either by manual or digital tracing) is assumed to be immediately diagnosed. This neglects a crucial difference between manual and digital tracing: manual tracing is time- and labor-intensive, so there is often a delay between the case confirmation and notification of contacts, whereas digital tracing is designed to avoid (or shorten) this delay. Additionally, 
a limitation of a contact tracing app is that there could be close contacts that can lead to transmission but are too short to be registered by the app. For this reason, it would be of interest to extend our model to allow for infectious contacts between app-users that are not registered with a certain (small) probability and thus not leading to tracing. It is however not straightforward to analyze such an extended model since it somehow mixes aspects of both digital (requiring both to be app-users) and manual (the possibility that such an app-contact is not recorded) contact tracing. Further, we make the simplifying assumption that traced individuals who have by then recovered are also identified and contact traced (see \cite{muller_contact_2000} which does not make this assumption). In reality, this may happen for individuals who have recently recovered, but perhaps not for those who have recovered several weeks earlier. Therefore, our results in this paper give an upper bound on how effective the real contact tracing would be. 

A different extension would be to consider a structured community (e.g. a social network) combined with homogeneous random contacts. Then, digital contact tracing would work in the same way for random contacts and contacts in the social network (as long as both are app-users), whereas manual contact tracing would most likely only happen on the social network (very rarely could you e.g. name whom you sat next to in the bus).

On the mathematical side, we only focus on the early stage of the epidemic in this paper. Analyzing later parts of the epidemic remains open. Moreover, obtaining a more explicit expression for $ R_{DM}$ in the combined model could shed more light on the effect of combining digital and manual contact tracing. For instance, it was shown numerically not to be monotone in $\pi$ and $p$, but one would expect the critical curve $R_{DM}=1$ to be monotone in the $(\pi ,p)$ plane (cf. Figure \ref{fig:RDM} and \ref{fig:RDMsmall}), but a proof of this missing. Further, it would be of interest to derive an individual reproduction number $R^{(ind)}_{DM}$ and prove properties for this reproduction number, e.g.\ that it has the same threshold as $R_{DM}$ and that it decreases monotonically in both $\pi$ and $p$.

Despite these model limitations, we consider our contribution an important step in increasing understanding of the effect of different contact tracing forms on spreading disease.

\acks We thank Frank Ball for his idea to derive a simplified expression of component reproduction number $R^{(c)}_{M}$ for manual tracing in \cite{zhang_analysing_2021}, which also made it possible to simplify the reproduction number $R_{D}$ in the current paper. 
\fund 
T.B. is grateful to the Swedish Research Council (grant 2020-04744) for financial support.

\competing 
There were no competing interests to declare which arose during the preparation or publication process of this article.



%
%
%

\end{document}